\documentclass[11pt]{article}
\usepackage{}
\usepackage{enumerate}
\usepackage{mathbbold}
\usepackage{bbm}
\usepackage{amsfonts}
\usepackage[T1]{fontenc}
\usepackage{latexsym,amssymb,amsmath,amsfonts,amsthm}
\usepackage{graphicx, color}
\usepackage{subfigure}
\usepackage{epsf}
\usepackage{epstopdf}
\usepackage{amssymb}
\usepackage{mathrsfs}

\renewcommand{\theequation}{\arabic{section}.\arabic{equation}}

\newcommand{\R}{{\mathbb R}}

\newcommand{\be}{\begin{eqnarray}}
\newcommand{\ben}{\begin{eqnarray*}}
\newcommand{\en}{\end{eqnarray}}
\newcommand{\enn}{\end{eqnarray*}}
\newcommand{\ba}{\backslash}

\newcommand{\ov}{\overline}

\newcommand{\G}{\Gamma}

\newcommand{\om}{\omega}

\newcommand{\hx}{\hat{x}}
\newtheorem{theorem}{Theorem}[section]
\newtheorem{lemma}[theorem]{Lemma}

\begin{document}
\title{\bf A multi-frequency sampling method for the inverse source problems with  sparse measurements}
\author{Xiaodong Liu\footnotemark[1] \quad and \quad Shixu Meng\footnotemark[2]}
\renewcommand{\thefootnote}{\fnsymbol{footnote}}
\footnotetext[1]{Academy of Mathematics and Systems Science, Chinese Academy of Sciences, Beijing 100190, China. {\tt xdliu@amt.ac.cn}}
\footnotetext[2]{Corresponding author. Academy of Mathematics and Systems Science, Chinese Academy of Sciences,
Beijing 100190, China.  {\tt shixumeng@amss.ac.cn}}
\maketitle

\begin{abstract}
We consider the  inverse source problems with multi-frequency  sparse near field measurements.
In contrast to the existing near field operator based on the integral over the space variable,  a multi-frequency near field operator  is introduced based on the integral over the frequency variable.  A factorization of this multi-frequency near field operator is further given and analysed. Motivated by such a factorization, we introduce a multi-frequency sampling method to reconstruct the source support. Its theoretical foundation is then derived from the properties of the factorized operators and a properly chosen point spread function. Numerical examples are provided to illustrate the multi-frequency sampling method with sparse near field measurements. Finally we briefly discuss how to extend the  near field case to the far field case.

\vspace{.2in}
{\bf Keywords:} sampling method; multi-frequency; sparse data; inverse source problems

\vspace{.2in} {\bf AMS subject classifications:}
35P25, 45Q05, 78A46, 74B05

\end{abstract}

\section{Introduction}
Inverse scattering merits numerous applications in radar, medical imaging, geophysics, non-destructive testing, and many others. In the last thirty years, sampling methods for shape reconstruction in inverse scattering problems have attracted a lot of interest.
Classical examples include the linear sampling method  \cite{ColtonKirsch} and the factorization method  \cite{Kirsch98}.  Such methods are independent of any a priori information on the geometry and physical properties of the unknown scatterers.  The theoretical foundation of these sampling methods  is the factorization of the far field or the near field operator $F$   in the form
\be\label{Ffacobstacles-inhomogeneity}
F = \tilde{H} TH
\en
with some coercive operator $T$ and a pair of  operators $H$ and $\tilde{H}$ (in some cases $\tilde{H}$ is the adjoint of $H$). The operators $H$, $\tilde{H}$, and $T$ may vary in different inverse scattering problems. With this factorization, one may derive a necessary and sufficient criteria to reconstruct the scatterer and design an indicator function which is large inside the underlying scatterer and relatively small outside.
We refer to the monographs  \cite{CaCo,CK,kirsch2008factorization} for a comprehensive introduction.

There have been recent efforts on other types of sampling methods  such as orthogonality sampling  method \cite{griesmaier2011multi,harris2020orthogonality,Potthast2010}, direct sampling method  \cite{ArensJiLiu,BaoHuangLiZhao,ItoJinZou,JiLiu-electromagnetic,LiZou}, single-shot method   \cite{LiLiuZou},
reverse time migration \cite{CCHuang}, and other direct reconstruction methods \cite{ammari2013mathematical,ammari2012direct}.
These sampling methods inherit many advantages of the classical ones. The main feature of these  sampling methods is that only inner product of the measurements with some suitably chosen functions is involved in the imaging function and thus  these  sampling methods are robust to noises. Recently,  a sampling type method is of the form \cite{LiuIP17,LiuMengZhang}
\begin{equation}
I(z) = \left| ( F g_z,\tilde{g}_z) \right|
\end{equation}
with suitably chosen $g_z$ and $\tilde{g}_z$ such that
\begin{itemize}
\item $\tilde{H}^* \tilde{g}_z=H g_z$ where $\tilde{H}^*$ is the adjoint of $\tilde{H}$ given in \eqref{Ffacobstacles-inhomogeneity};
\item $|Hg_z|$ is a point spread function   that peaks when the sampling point $z$ located on the boundaries of  the unknown scatterers (or inside the scatterers in different problems).
\end{itemize}
Based on the factorization of the far field or near field operator \eqref{Ffacobstacles-inhomogeneity},
we deduce that their exists two positive constant $c$ and $C$ such that
\be \label{intro foundation}
c\|Hg_z\|^2\leq I(z)\leq C \|Hg_z\|^2.
\en
This theoretical foundation \eqref{intro foundation} implies that $I(z)$ is qualitatively the same as a suitable norm of the point spread function $|Hg_z|$, i.e.
\ben
I(z)\approx \|Hg_z\|^2,
\enn
which allows us to design a robust and fast imaging algorithm.

In this work, we apply the above framework to introduce a multi-frequency sampling method for the inverse source problems  using  sparse measurements.
Most of the works in literature assume that the measurements are taken all around the unknown source. In this full-aperture case,  uniqueness of the source function $f$ can be proved using the Fourier transform, see for instance \cite{BaoLinTriki, EllerValdivia}. We refer to \cite{BaoLuRundellXu} for an iterative method with respect to frequencies.
Recently, some non-iterative methods have been proposed. In particular, \cite{WangGuoZhangLiu2017IP, ZhangGuo}  studied a Fourier method to reconstruct the source function and \cite{BousbaGuoWangLi,ZGLL} investigated the direct sampling methods for identifying the point sources.
One practical situation involves sparse measurements, i.e. measurements that are only  available at a limited number of sensors. Generally, neither the source nor its support is uniquely determined by the sparse data. Some lower bounds for the support of a source have been established in \cite{SylvesterKelly} in terms of multi-frequency sparse far field measurements. If it is known a priori that the source is a combination of point sources,  uniqueness for the locations and the scattering strengths can be established \cite{Ji-multipolarFarfield, JiLiu-point, JiLiu-nearfield}. We refer  to \cite{GriesmaierSchmiedecke-source} and \cite{AlaHuLiuSun} for a factorization method and a direct sampling method, respectively, using  multi-frequency sparse far field measurements. For the multi-frequency sparse near field measurements, the corresponding uniqueness  and direct sampling method can be found in a recent work \cite{JiLiu-nearfield}.
Following the framework \eqref{Ffacobstacles-inhomogeneity}--\eqref{intro foundation} developed in \cite{LiuIP17,LiuMengZhang} for inverse obstacle/medium scattering problems, this work contributes to a multi-frequency sampling method for the inverse source problems with  sparse measurements.

The paper is further organized as follows. We introduce the model of the forward problem and inverse source problem with multi-frequency sparse near field measurements in Section \ref{section model}. In contrast to the existing near field operator based on the integral over the space variable, we define a multi-frequency near field operator for each measurement point based on the integral over the frequency variable. We further provide a factorization of this multi-frequency near field operator in the form \eqref{Ffacobstacles-inhomogeneity} in Section \ref{section factorization}. This factorization is one of the key ingredients in our  multi-frequency sampling method.
In Section \ref{section imaging}, we design an imaging function with  properly chosen $g_z$ and $\tilde{g}_z$. Furthermore we provide its theoretical foundation based on the factorization and a point spread function. Section \ref{section numerical examples} presents numerical examples to illustrate our multi-frequency sampling method with sparse near field measurements. Finally in Section \ref{section appendix} we briefly discuss how to extend the  near field case to the far field case.

\section{Mathematical model and the inverse problem} \label{section model}
We consider the acoustic wave propagation due to a source $f$ in a homogeneous isotropic medium in $\R^3$ with speed of sound $c>0$. Denote by $D$ the support of the unknown source, which is a bounded Lipschitz domain in $\R^3$ with connected complement $\R^3\ba\ov{D}$.
The mathematical model of the scattering of time-harmonic wave leads to the nonhomogeneous Helmholtz equation
\be\label{Helmholtzequation}
\Delta u^s +k^2 u^s = f \quad \mbox{in}\,\, \R^3,
\en
where the  wave number is $k=\om/c$ with frequency $\om>0$.
The scattered field $u^s$ is required to satisfy the Sommerfeld radiation condition
\be\label{SRC}
\lim_{r=|x|\rightarrow \infty} r \Big(\frac{\partial u^s}{\partial r} - iku^s \Big) = 0
\en
uniformly for all directions $x/|x|$.
Throughout the paper we write $u^s(\cdot, k)$ to emphasize the dependence on the wave number $k$.
If not otherwise stated, we consider multiple wave numbers in a bounded interval, i.e.,
\be\label{kassumption}
k\in K:=(0, k_{max}],
\en
with an upper bound $k_{max}>0$.


The unique radiating solution $u^s$ to the scattering problem \eqref{Helmholtzequation}-\eqref{SRC} takes the form
\be\label{us}
u^s(x,k)=\int_{D}\Phi_{k}(x,y)f(y)ds(y),\quad x\in\R^3, \, k \in K
\en
with
\ben \label{FundaSolu_of_Helmholtz}
\Phi_{k}(x,y):= \frac{e^{ik|x-y|}}{4\pi|x-y|}, \quad x\not=y
\enn
being the fundamental solution to the Helmholtz equation.

A practical inverse problem is to recover the source support $D$ from the multi-frequency sparse measurements where the scattered fields are available at finitely many points in the set
\ben
\G_L:=\{x_1, x_2, \cdots, x_L\}\subset \R^3\ba\ov{D}.
\enn
Accordingly, we obtain the following multi-frequency sparse near field measurements
\be\label{eqivalentdatanear}
\mathbb{M}_{N}:=\{u^s(x, k)\,|\, x\in\G_L,\,\, k\in K\}.
\en

Generally speaking, uniqueness of the source support $D$ does not hold from the near field measurements $\mathbb{M}_{N}$. An approximate  support   containing $D$ can be uniquely determined. We refer to \cite{AlaHuLiuSun,GriesmaierSchmiedecke-source,JiLiu-nearfield,SylvesterKelly} for more details on the uniqueness.


\section{The multi-frequency near field operator and its factorization} \label{section factorization}
In this section, we introduce a multi-frequency near field operator given by the measurements $\mathbb{M}_N$ \eqref{eqivalentdatanear} and study its factorization which plays a key role in the analysis and design of our imaging algorithm.
We assume that the source function $f\in L^2(\R^3)$ is real-valued and satisfies that
\be\label{Assumpf}
C_f\geq f(x)\geq c_f \quad \mbox{or}\quad -C_f\leq f(x) \leq -c, \quad x\in\, D
\en
for two constants $C_f > c_f >0$.

To begin with, we define the scattered fields with negative frequencies by
\ben
u^s(x, -k):= \ov{u^s(x, k)}, \quad x\in \G_L, \, k\in K.
\enn
Note that $f$ is real-valued, we then have the following integral representation for the scattered fields
\be\label{us complex}
u^s(x,k)=\int_{D}\Phi_{k}(x,y)f(y)ds(y),\quad x\in \G_L, \, k \in  [-k_{max}, k_{max}]\ba\{0\}.
\en

For each $x\in\G_L$, we define the integral operator $\mathcal {N}_{x}: L^2(K)\rightarrow L^2(K)$ by
\be\label{NearFieldOperator}
(\mathcal {N}_{x}g)(t) := \int_{K} u^s(x, t-s)g(s)ds, \quad t\in K,
\en
which we will call the {\em multi-frequency near field operator}. Different from the usual near field operator which uses full-aperture measurements with a single frequency, this operator seeks to explore the relations among different  frequencies. The factorization of the usual near field operator plays an important role in the modified sampling method \cite{LiuMengZhang} so that no asymptotic assumptions on the distance between the measurement surface and the scatterers were made there. In the same spirit, we seek to given a characterization of the unknown source support $D$ by the multi-frequency near field operator $\mathcal {N}_{x}$ together with its factorization.


To factorize $\mathcal{N}_x$, we first define $\mathcal {P}_{x}: L^2(D)\rightarrow L^2(K)$ by
\be \label{P}
(\mathcal {P}_{x} \psi)(t) := \int_{D} e^{it|x-y|}\psi(y) dy, \quad t\in K.
\en
We can directly derive that its adjoint $\mathcal {P}^{\ast}_{x}: L^2(K)\rightarrow L^2(D)$ is given by
\be \label{P*}
(\mathcal {P}^{\ast}_{x} \phi)(y) := \int_{K} e^{-is|x-y|}\phi(s) ds, \quad y\in D.
\en
Finally we define the operator $\mathcal{T}_{x}: L^2(D)\rightarrow L^2(D)$ by
\be\label{Tx}
(\mathcal {T}_{x}h)(y):= \frac{f(y)}{4\pi|x-y|}h(y), \quad y\in D.
\en
It directly follows that $\mathcal {T}_x$ is a self-adjoint operator since $f$ is real-valued.
Now we are ready to give the following factorization of the multi-frequency near field operator $\mathcal{N}_x$.

\begin{theorem}\label{Thm-Nfac}
For each $x \in \G_L$, the following factorization holds
\be\label{Nfactorization}
\mathcal {N}_x = \mathcal {P}_{x} \mathcal {T}_x \mathcal {P}^{\ast}_{x}.
\en
\end{theorem}
\begin{proof}
Inserting the integral representation \eqref{us complex} into the definition \eqref{NearFieldOperator} of the multi-frequency  near field operator $\mathcal {N}_{x}$, we have that
\ben
(\mathcal {N}_{x}g)(t)
&=& \int_{K} \int_{D} \frac{e^{i(t-s)|x-y|}}{4\pi|x-y|}f(y)dy g(s)ds\cr
&=& \int_{D} e^{it|x-y|}\frac{f(y)}{4\pi|x-y|}\int_{K}e^{-is|x-y|} g(s)dsdy\cr
&=& (\mathcal {P}_{x} \mathcal {T}_{x} \mathcal {P}^{\ast}_{x}g)(t),\quad t\in K.
\enn
This completes the proof.
\end{proof}
Theorem \ref{Thm-Nfac} provides a symmetric factorization which allows us to study a  factorization method, which is expected to be similar to the factorization method using multi-frequency far field measurements \cite{GriesmaierSchmiedecke-source}. Alternatively, our goal here is to utilize this factorization to introduce a multi-frequency sampling method. We remark that this type of sampling method may still work even in the case of non-symmetric factorization (see for instance \cite{LiuMengZhang}).

%
\section{The  multi-frequency sampling method} \label{section imaging}
This section is devoted to  a  multi-frequency sampling method using sparse near field measurements $\mathbb{M}_{N}$ \eqref{us}. To begin with, we provide the following property of the middle operator $\mathcal{T}_x$.

\begin{lemma}\label{Itheorey lemma}
Assume that $f\in L^2(\R^3)$ satisfies \eqref{Assumpf}. For each $x \in \G_L$, the operator $\mathcal {T}_x: L^2(D)\rightarrow L^2(D)$ is  self-adjoint and coercive, i.e.
\be \label{Itheorey lemma eqn1}
\frac{c_f}{4\pi r_2}\|h\|^2_{L^{2}(D)} \leq |(\mathcal {T}_x h, h)_{L^2(D)}| \leq \frac{C_f}{4\pi r_1}\|h\|^2_{L^{2}(D)}, \quad \forall h\in L^2(D),
\en
where
\be\label{r1r2}
r_1:=\inf\{|x-y|\, :\, y\in D\} \quad\mbox{and}\quad r_2:=\sup\{|x-y|\, :\, y\in D\}.
\en
\end{lemma}
\begin{proof}
Recall the definition of $\mathcal {T}_x$  \eqref{Tx}, we have that
$$
 (\mathcal {T}_x h, h)_{L^2(D)} = \int_D \frac{f(y)}{4\pi|x-y|} |h(y)|^2 dy.
$$
Note from the assumption \eqref{Assumpf}  that
\be
C_f\geq f(x)\geq c_f \quad \mbox{or}\quad -C_f\leq f(x) \leq -c, \quad x\in\, D,
\en
then the inequality \eqref{Itheorey lemma eqn1} follows by the definition of $r_1$ and $r_2$. This proves the lemma.
\end{proof}
We now prove the following theorem which plays an important role in the analysis and design of our imaging algorithm.
\begin{theorem}\label{Itheorey}
Assume that $f\in L^2(\R^3)$ satisfies \eqref{Assumpf}. For each $x \in \G_L$, the following inequality holds
\ben
\frac{c_f}{4\pi r_2}\|\mathcal {P}^{\ast}_{x}g\|^2_{L^{2}(D)} \leq |(\mathcal {N}_x g, g)_{L^2(K)}| \leq \frac{C_f}{4\pi r_1}\|\mathcal {P}^{\ast}_{x}g\|^2_{L^{2}(D)}, \quad \forall g\in L^2(K).
\enn
\end{theorem}
\begin{proof}
With the help of the factorization \eqref{Nfactorization} in Theorem \ref{Thm-Nfac}, we have that
\ben
(\mathcal {N}_x g, g)_{L^2(K)}= (\mathcal {P}_{x} T_x \mathcal {P}^{\ast}_{x} g, g)_{L^2(K)}= (T_x \mathcal {P}^{\ast}_{x} g, \mathcal {P}^{\ast}_{x}g)_{L^2(D)}.
\enn
Then the theorem follows from Lemma \ref{Itheorey lemma}. This proves the theorem.
\end{proof}

\begin{figure}
\begin{center}
\resizebox{1\textwidth}{!}{\includegraphics{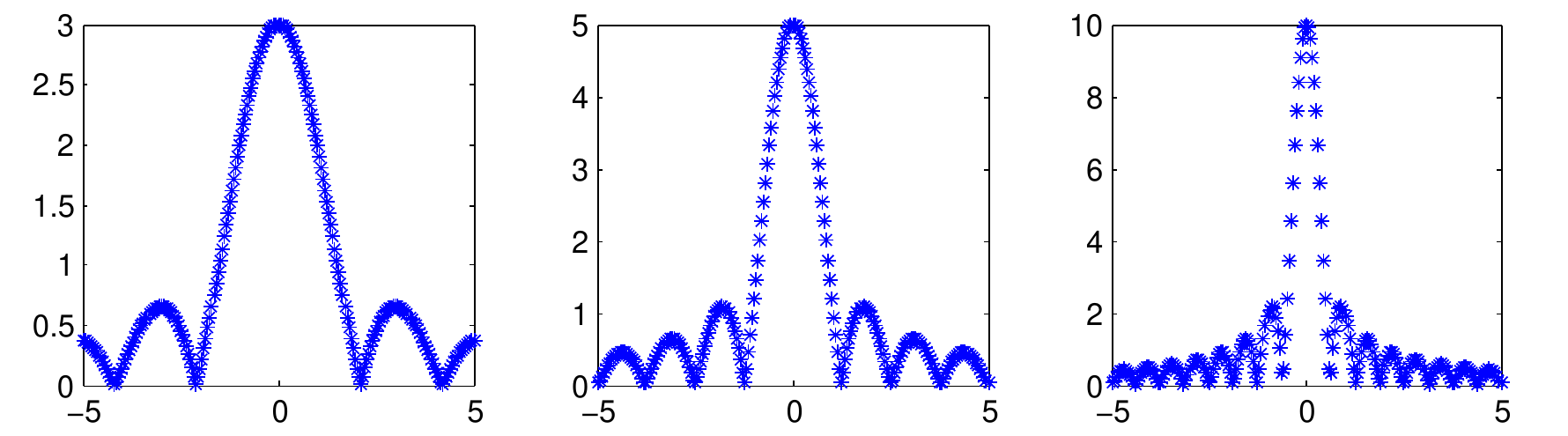}}
\end{center}
\caption{Plot of $|(\mathcal {P}^{\ast}_{x}g_{xz})(y)|$ as a function of $\mathbbm{t}$  in $[-5, 5]$. Left: $k_{max}=3$; \quad Middle: $k_{max}=5$; \quad Right: $k_{max}=10$.}
\label{fBehaviour}
\end{figure}

The above theorem implies that $|(\mathcal {N}_x g, g)_{L^2(K)}|$  is qualitatively the same as $\|\mathcal {P}^{\ast}_{x}g\|^2_{L^{2}(D)}$. The following theorem gives a properly chosen $g_{xz}$ such that $(\mathcal {P}^{\ast}_{x}g_{xz})$ is a point spread function.
\begin{theorem}
For each $x \in \G_L$ and  $z\in\R^3$, we define $g_{xz} \in L^2(K)$ by
\be\label{Def-gz}
g_{xz}(s):=e^{is|x-z|}, \quad s\in K.
\en
It holds that
\be \label{P* gz}
(\mathcal {P}^{\ast}_{x}g_{xz})(y)=
\Bigg\{
\begin{array}{cc}
 k_{max} &     \mbox{if }   \mathbbm{t} =0 \\
\frac{e^{i \mathbbm{t} k_{max} }-1}{\mathbbm{t} }  &  \mbox{if } \mathbbm{t}  \not=0
\end{array}
, \quad y\in D,
\en
where $\mathbbm{t} := |x-z|-|x-y|$. In particular, $(\mathcal {P}^{\ast}_{x}g_{xz})(y)=(\mathcal {P}^{\ast}_{x}g_{xz})(y^*)$ if $|y-x|=|y^*-x|$.

Furthermore,
$$
|(\mathcal {P}^{\ast}_{x}g_{xz})(y)| \le k_{max}, \quad y \in D,
$$
where ``$=$'' holds if and only if  $\mathbbm{t}=0$. The following asymptotic holds
$$
|(\mathcal {P}^{\ast}_{x}g_{xz})(y)| =\mathcal{O}\left(\frac{1}{\mathbbm{t}}\right), \quad \mathbbm{t} \rightarrow \infty.
$$
\end{theorem}
\begin{proof}
From the definition of $\mathcal {P}^{\ast}_{x}$ \eqref{P*}, we have that
\ben
(\mathcal {P}^{\ast}_{x}g_{xz})(y)=\int_{0}^{k_{max}}e^{is(|x-z|-|x-y|)}ds, \quad y\in D.
\enn
Then equation \eqref{P* gz} follows from a direct computation. It then   follows that $(\mathcal {P}^{\ast}_{x}g_{xz})(y)=(\mathcal {P}^{\ast}_{x}g_{xz})(y^*)$ if $|y-x|=|y^*-x|$.

Furthermore we derive that
\ben
\left|(\mathcal {P}^{\ast}_{x}g_{xz})(y)\right| = \left|\int_{0}^{k_{max}}e^{is\mathbbm{t}}ds\right| \le \int_{0}^{k_{max}}\left|e^{is\mathbbm{t}}\right| ds \le k_{max}, \quad y\in D,
\enn
where ``$=$'' holds if and only if  $\mathbbm{t}=0$.

Finally the asymptotic of $|(\mathcal {P}^{\ast}_{x}g_{xz})(y)|$ for large $\mathbbm{t}$   follows from \eqref{P* gz}. This proves the theorem.
\end{proof}


\begin{figure}
\begin{center}
\resizebox{0.3\textwidth}{!}{\includegraphics{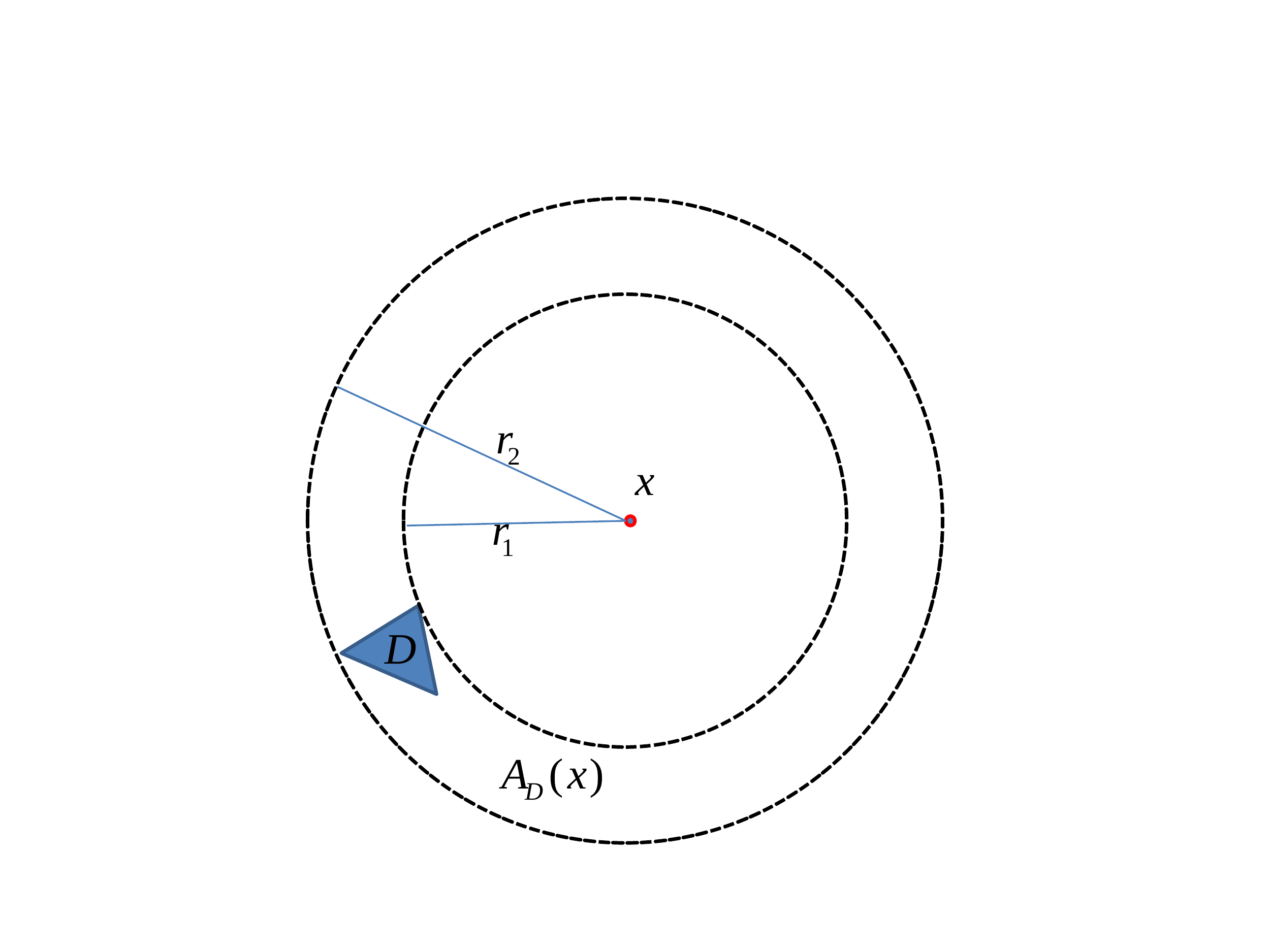}}
\end{center}
\caption{The smallest annulus support $A_D(x):=B_{r_2}(x)\ba\ov{B_{r_1}(x)}$ centered at the measurement point $x$ with $r_1=\inf\{|x-y|\, :\, y\in D\}$ and $r_2=\sup\{|x-y|\, :\, y\in D\}$.}
\label{SDAD}
\end{figure}
To illustrate the asymptotic property of $|(\mathcal {P}^{\ast}_{x}g_{xz})(y)|$ as a function of $\mathbbm{t}$, we refer to Figure \ref{fBehaviour} where we plot  $|(\mathcal {P}^{\ast}_{x}g_{xz})(y)|$  with respect to different $k_{max}$.  $|(\mathcal {P}^{\ast}_{x}g_{xz})(y)|$ is a point spread function of $\mathbbm{t}$.  In particular, Figure \ref{fBehaviour} implies that the resolution seems better with larger wave number $k_{max}$.

Now for   given $y \in D$ and sampling point $z$, $|(\mathcal {P}^{\ast}_{x}g_{xz})(y)|$ peaks at $\mathbbm{t}=0$, i.e. at sampling point $z$ such that $|x-z|=|x-y|$. In another word,  $|(\mathcal {P}^{\ast}_{x}g_{xz})(y)|$ as a function of sampling point $z$  peaks at a sphere centered at $x$ with radius $|x-y|$. Furthermore, turning to $\|\mathcal {P}^{\ast}_{x}g_{xz}\|^2_{L^{2}(D)}$ as a  function of sampling point $z$, we expect the function $\|\mathcal {P}^{\ast}_{x}g_{xz}\|^2_{L^{2}(D)}$ peaks in the smallest annulus support $A_D(x):=B_{r_2}(x)\ba\ov{B_{r_1}(x)}$  centered at the measurement point $x$ with $r_1$ and $r_2$ given by \eqref{r1r2}. See Figure \ref{SDAD} for an illustration.

Unfortunately, it is not feasible to compute $\|\mathcal {P}^{\ast}_{x}g_{xz}\|^2_{L^{2}(D)}$ directly since $D$ is what we aim to reconstruct. However we are able to relate $\|\mathcal {P}^{\ast}_{x}g_{xz}\|^2_{L^{2}(D)}$ to an indicator function given by measurements. Precisely, letting  $g=g_{xz}$ in Theorem \ref{Itheorey} yields that
\ben
\frac{c_f}{4\pi r_2}\|\mathcal {P}^{\ast}_{x}g_{xz}\|^2_{L^{2}(D)} \leq |(\mathcal {N}_x g_{xz}, g_{xz})_{L^2(K)}| \leq \frac{C_f}{4\pi r_1}\|\mathcal {P}^{\ast}_{x}g_{xz}\|^2_{L^{2}(D)},
\enn
i.e. $\|\mathcal {P}^{\ast}_{x}g_{xz}\|^2_{L^{2}(D)}$ is qualitatively the same as $|(\mathcal {N}_x g_{xz}, g_{xz})_{L^2(K)}|$ for a single measurement point $x$. This implies that  the measurement driven function $|(\mathcal {N}_x g_{xz}, g_{xz})_{L^2(K)}|$ peaks in the annulus support $A_D(x)$ without the knowledge of $D$.
Consequently for all $x\in \G_L$, we introduce the indicator function
\be\label{I-Near}
I(z):=\sum_{x\in\G_L}\Big|(\mathcal {N}_x g_{xz}, g_{xz})_{L^2(K)}\Big|,\quad z\in\R^3.
\en
To summarize, the indicator function $I(z)$ peaks  in the annulus support $A_D(x)$ for a single measurement point $x$. With the increase of   measurement points, an approximate source support $\bigcap_{x\in\G_L}A_D(x)$ is expected to be reconstructed by plotting $I(z)$ over a sampling region.


\section{Numerical examples} \label{section numerical examples}
In this section, we present some numerical examples to illustrate the performance of the multi-frequency sampling method.

\begin{figure}[htbp]
  \centering
    \includegraphics[width=8cm]{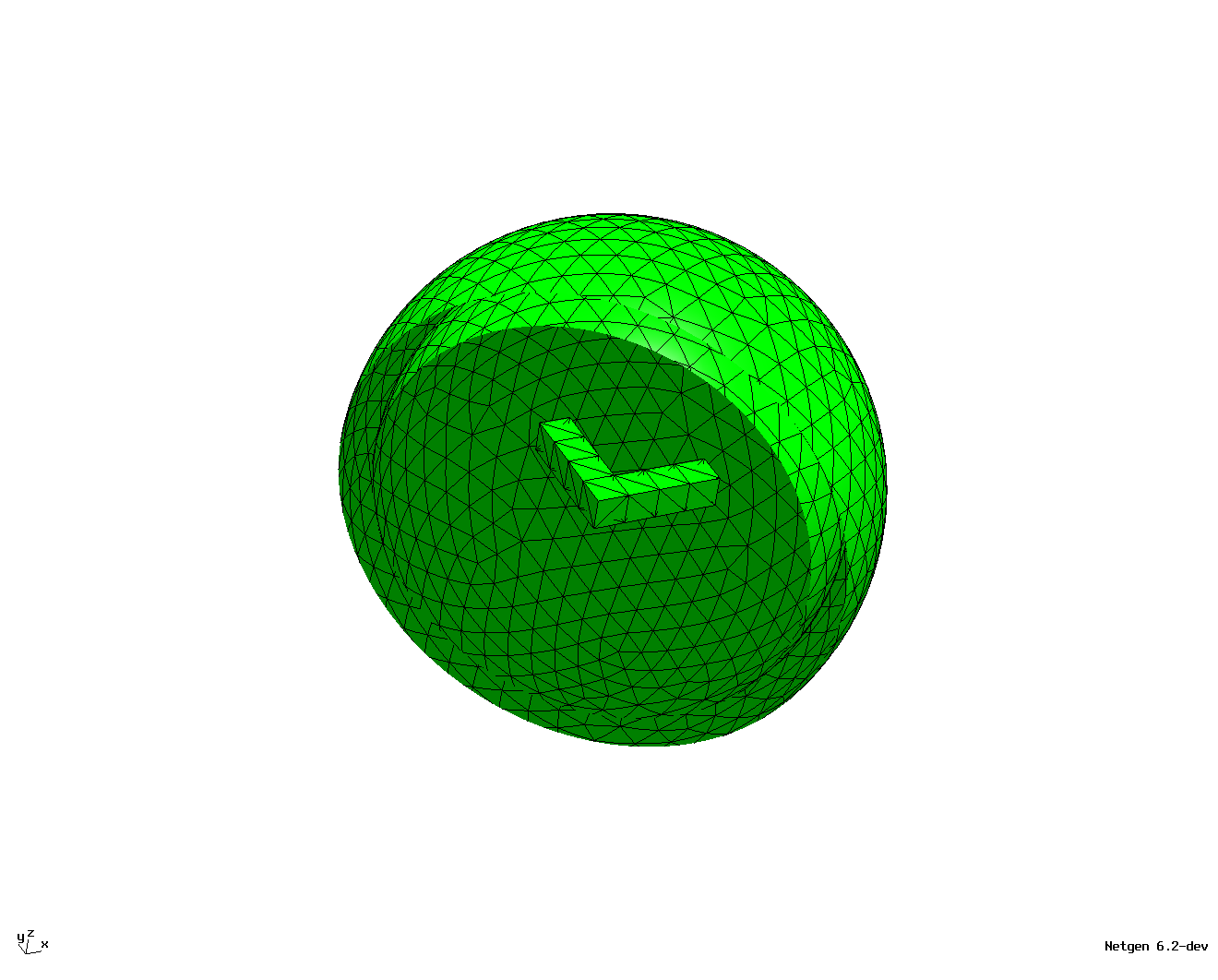}
\caption{{\bf Mesh view of the computational domain.} }
\label{domain}
\end{figure}

We generate the synthetic data $u^s$ using the finite element computational software Netgen/NGSolve \cite{schoberl1997netgen}. To be more precise, the computational domain is $\{x: |x|<4\}$ and the measurements are on the sphere $ \{x: |x|=3\}$. We apply a radial Perfectly Matched Layer (PML) in the domain $\{x: 3.5<|x|<4\}$ and choose PML absorbing coefficient $5i$. See Figure \ref{domain} for an illustration. In all of the numerical examples, we apply the second order finite element to solve for the  wave field where the source is constant $1$; the mesh size is chosen as $0.5$; the set of wave numbers is $\{k: k =1,2,3,\cdots,11\}$.  We further add $5\%$ Gaussian noise to the synthetic data $u^s$ to implement the indicator function given by \eqref{I-Near} in Matlab; for the visualization, we plot $I(z)$ over the sampling region $\{x=(x_1,x_2,x_3): |x_j|<3, \, j=1,2,3\}$ and we always normalize it such that its maximum value is $1$.

\begin{figure}[htbp]
  \centering
  \includegraphics[width=5cm]{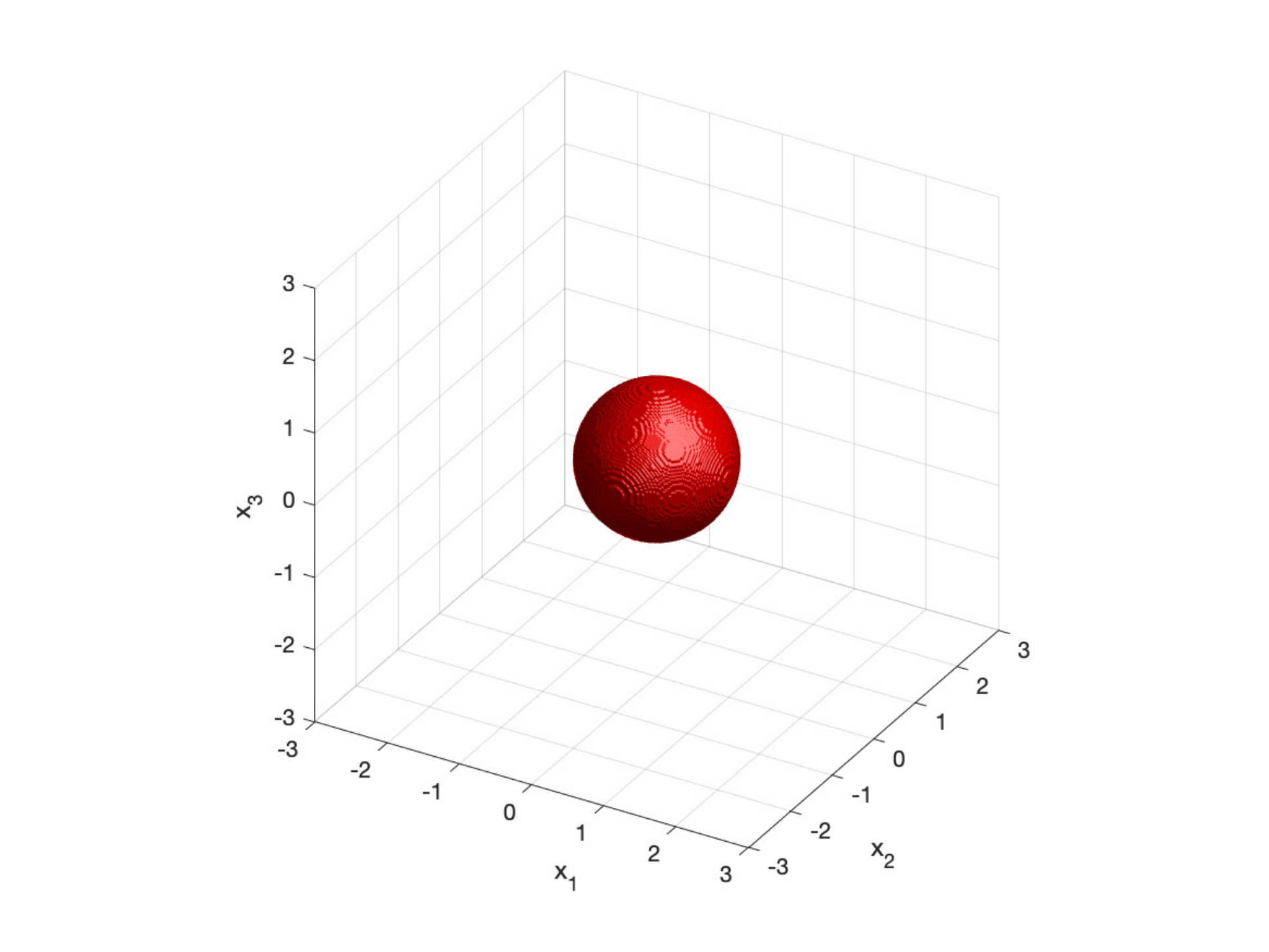}
    \includegraphics[width=5cm]{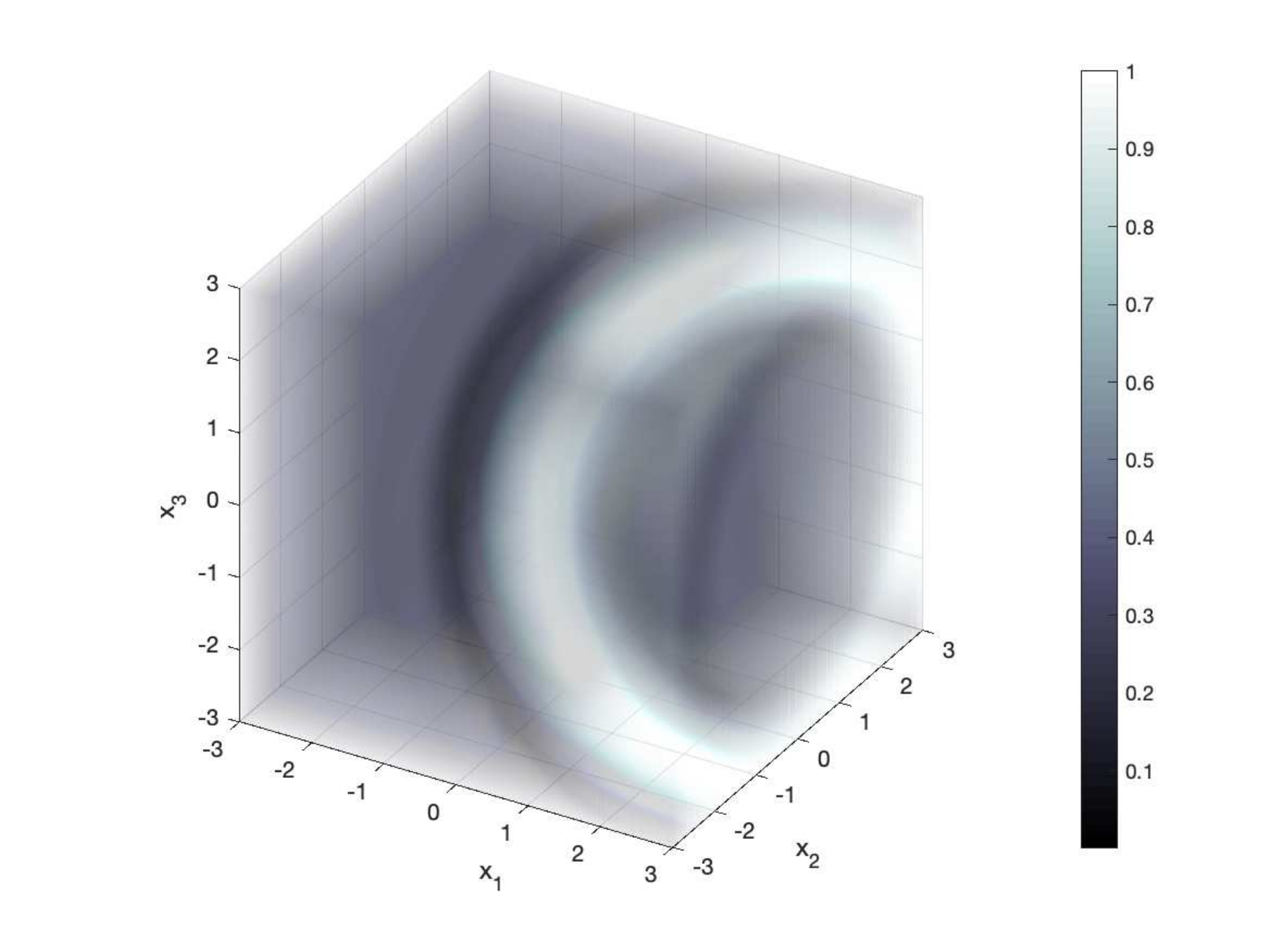}
     \includegraphics[width=5cm]{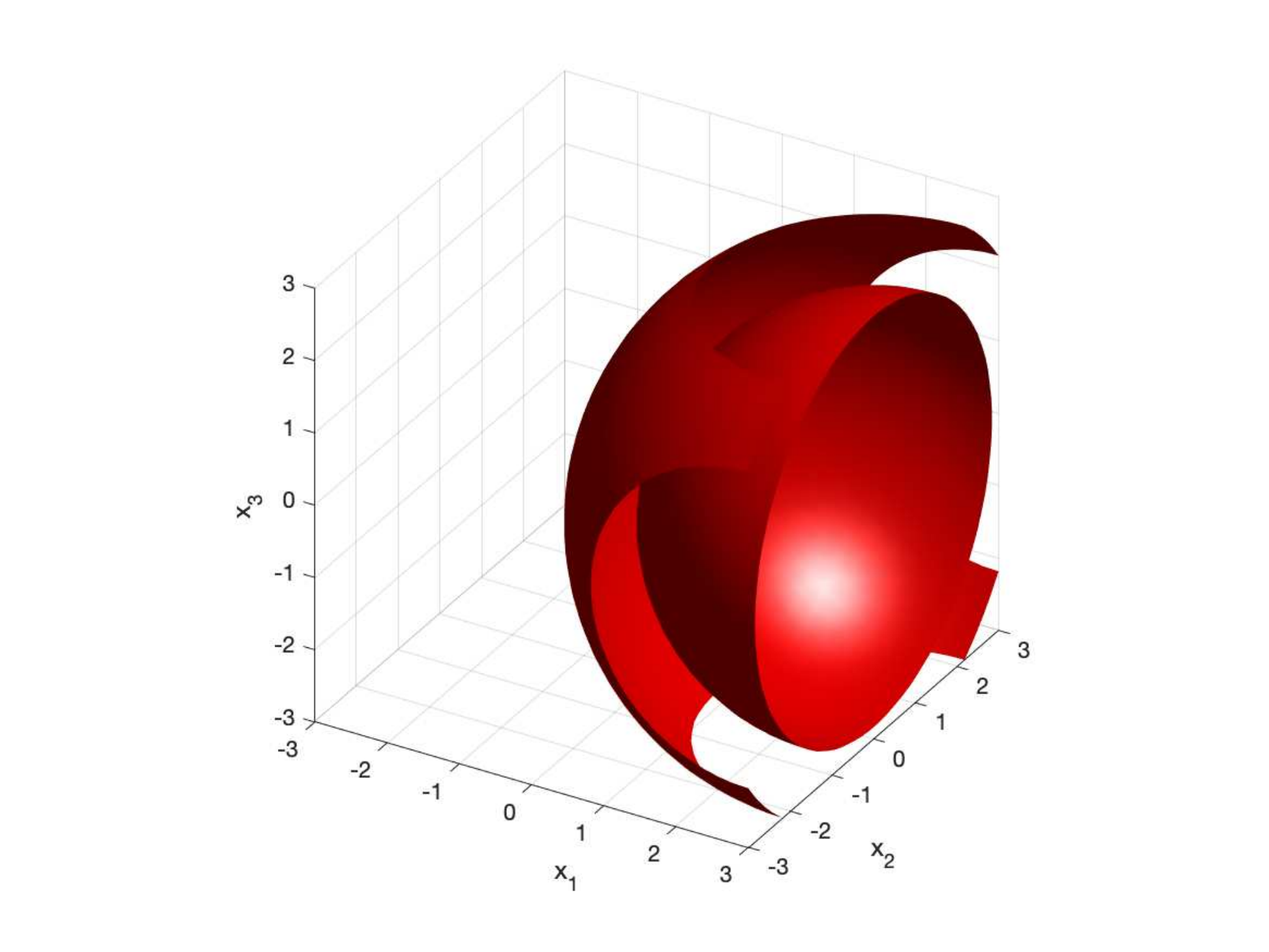}
\caption{{\bf Reconstructions of a ball with one measurement point $(3,0,0)$.} Left: 3D view of the exact ball. Middle: 3D view of the reconstruction. Right: iso-surface view of the reconstruction with iso-value $7\times 10^{-1}$.}
\label{ball one point}
\end{figure}

We first illustrate the performance of the multi-frequency sampling method with one measurement point. In this example, the measurement point is   $(3,0,0)$. The source has support given by a ball
\begin{eqnarray} \label{numeric ball def}
\{x: |x_1|^2 +  |x_2|^2 +  |x_3|^2< 1 \}.
\end{eqnarray}
We plot the exact ball,  the three dimensional view of the reconstruction, and its iso-surface view with iso-value $7\times 10^{-1}$ in Figure  \ref{ball one point}. As expected, the reconstruction is an annulus support of the exact ball.

\begin{figure}[htbp]
  \centering
  \includegraphics[width=5cm]{figures/Ball_exact-eps-converted-to.pdf}
  \includegraphics[width=5cm]{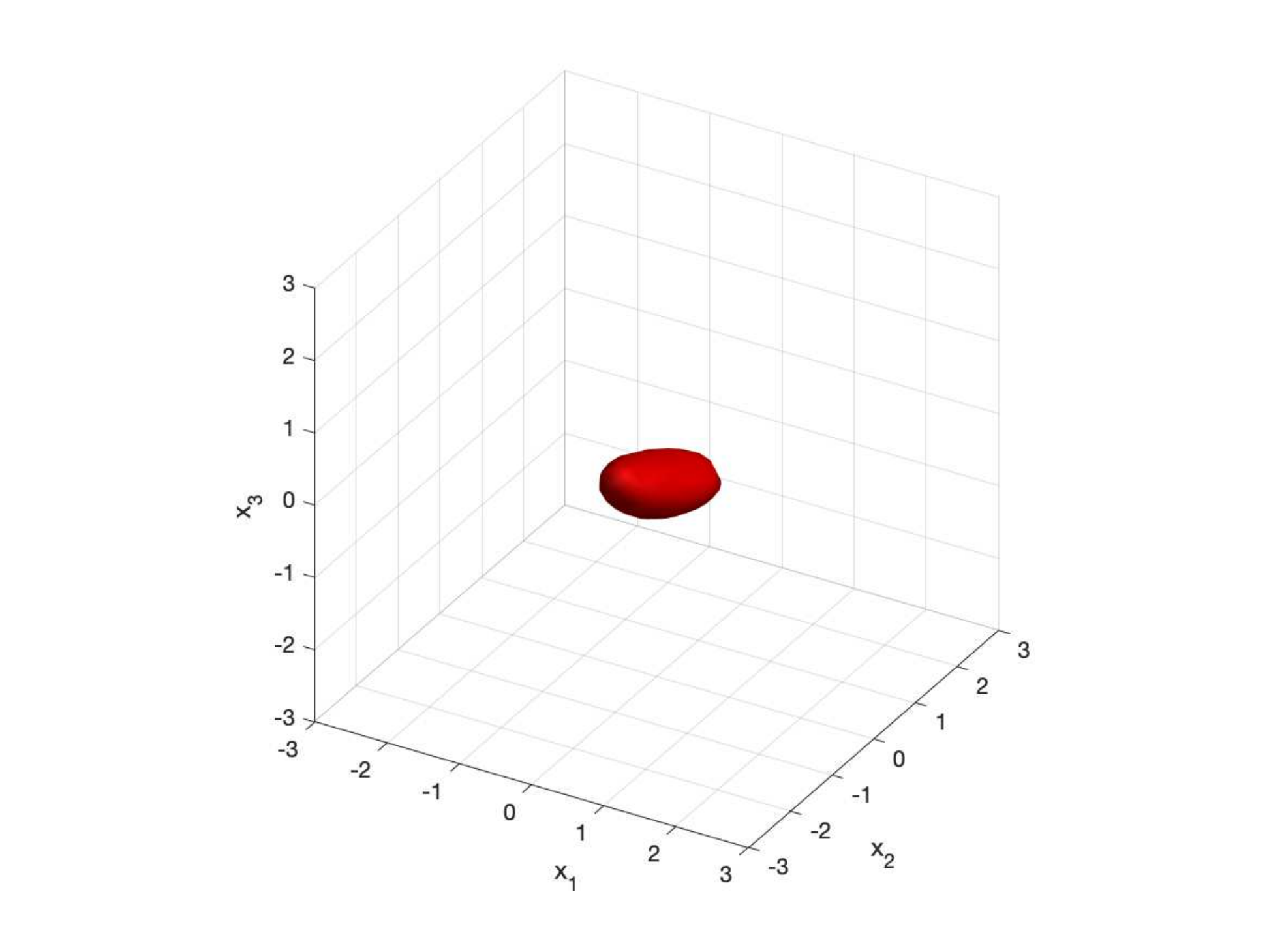}
    \includegraphics[width=5cm]{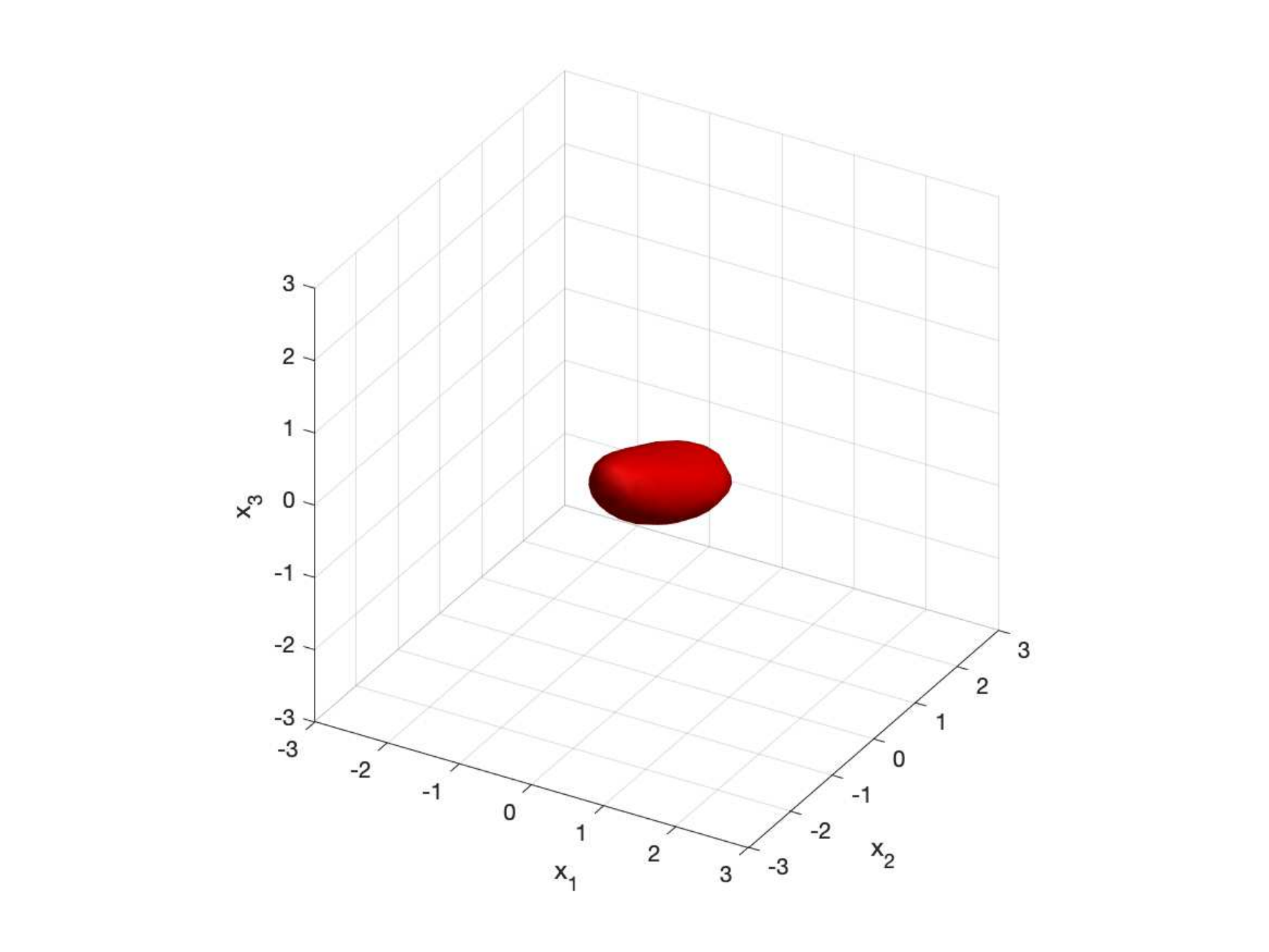}
    \includegraphics[width=5cm]{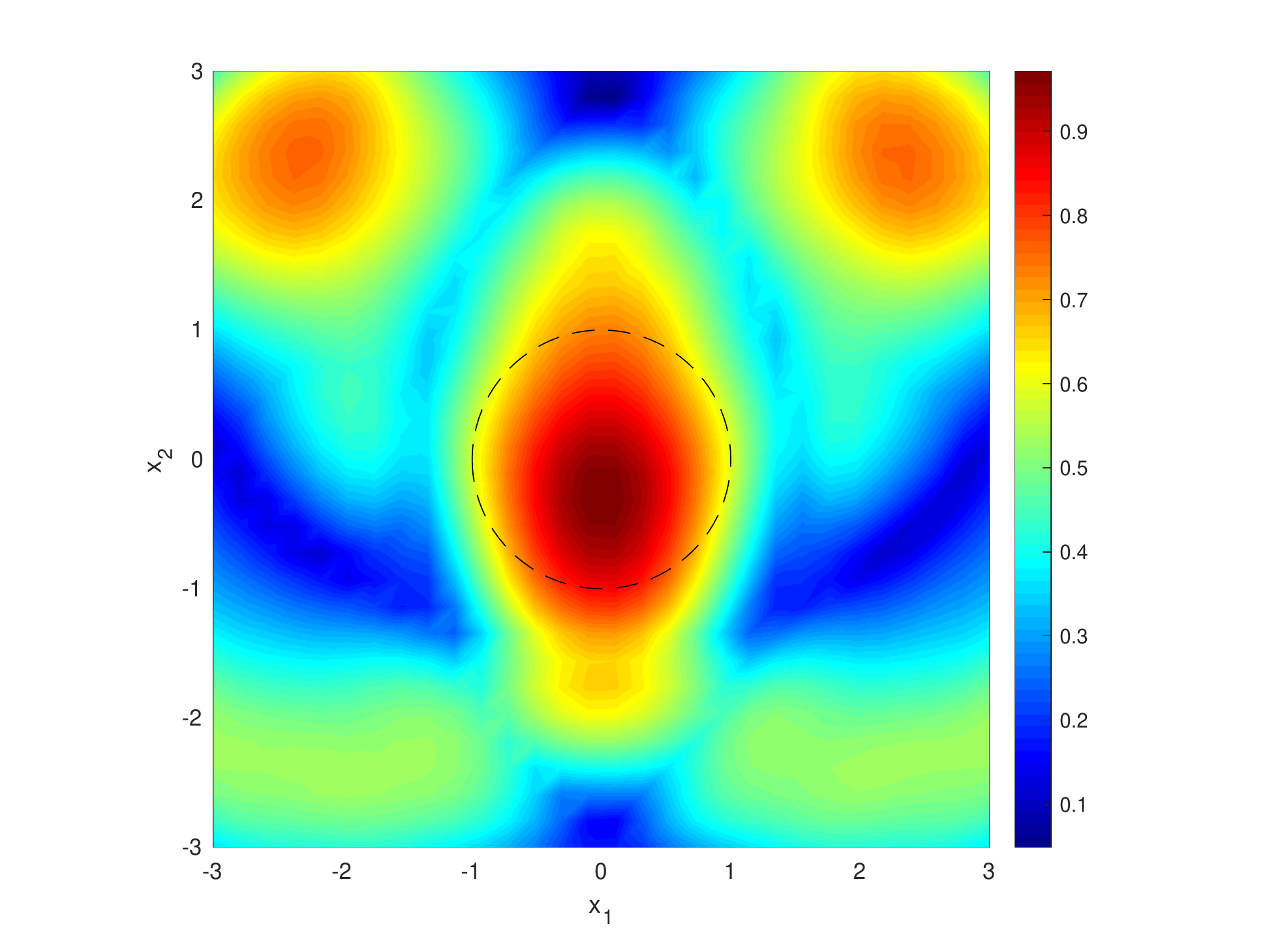}
     \includegraphics[width=5cm]{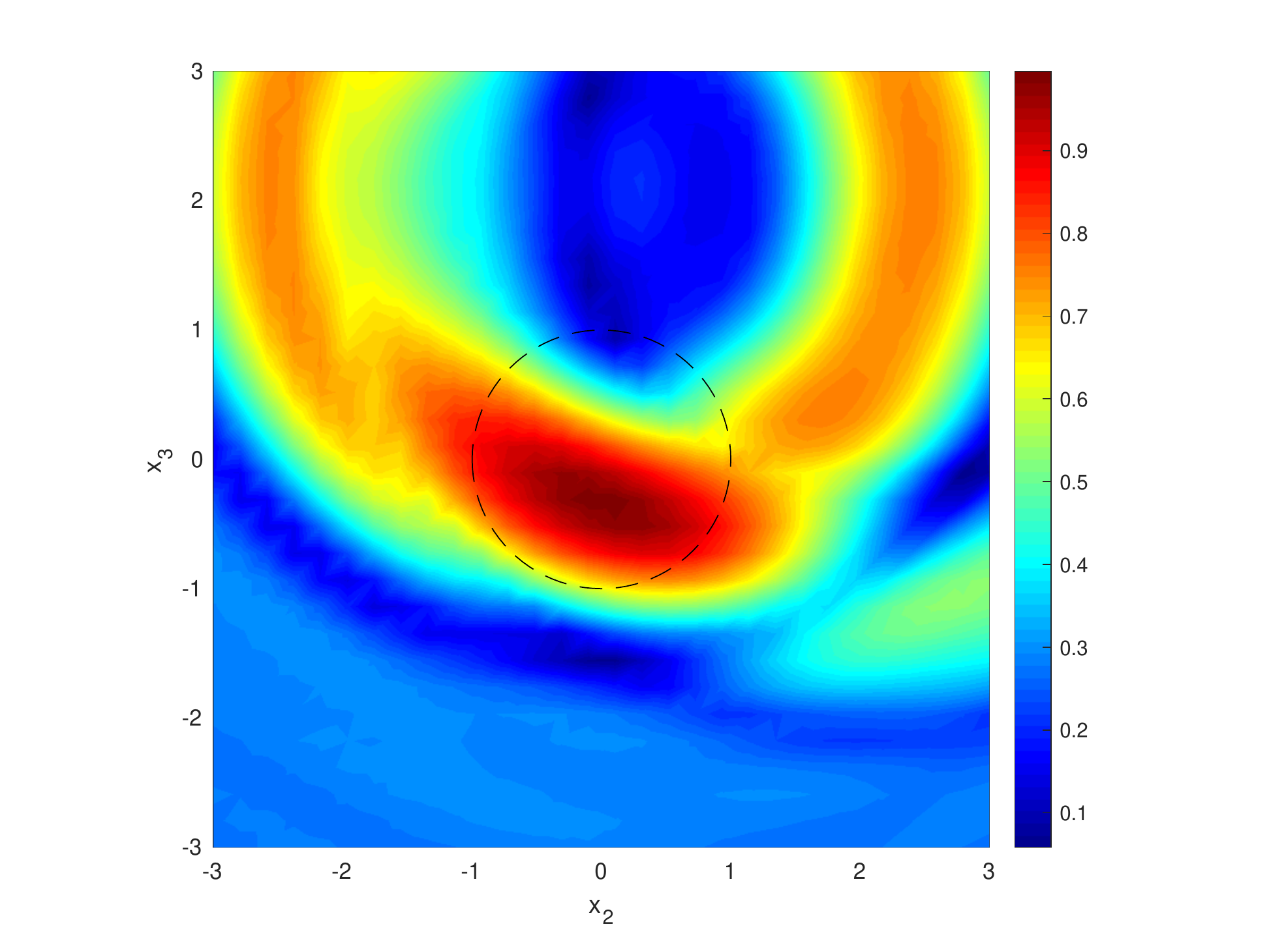}
          \includegraphics[width=5cm]{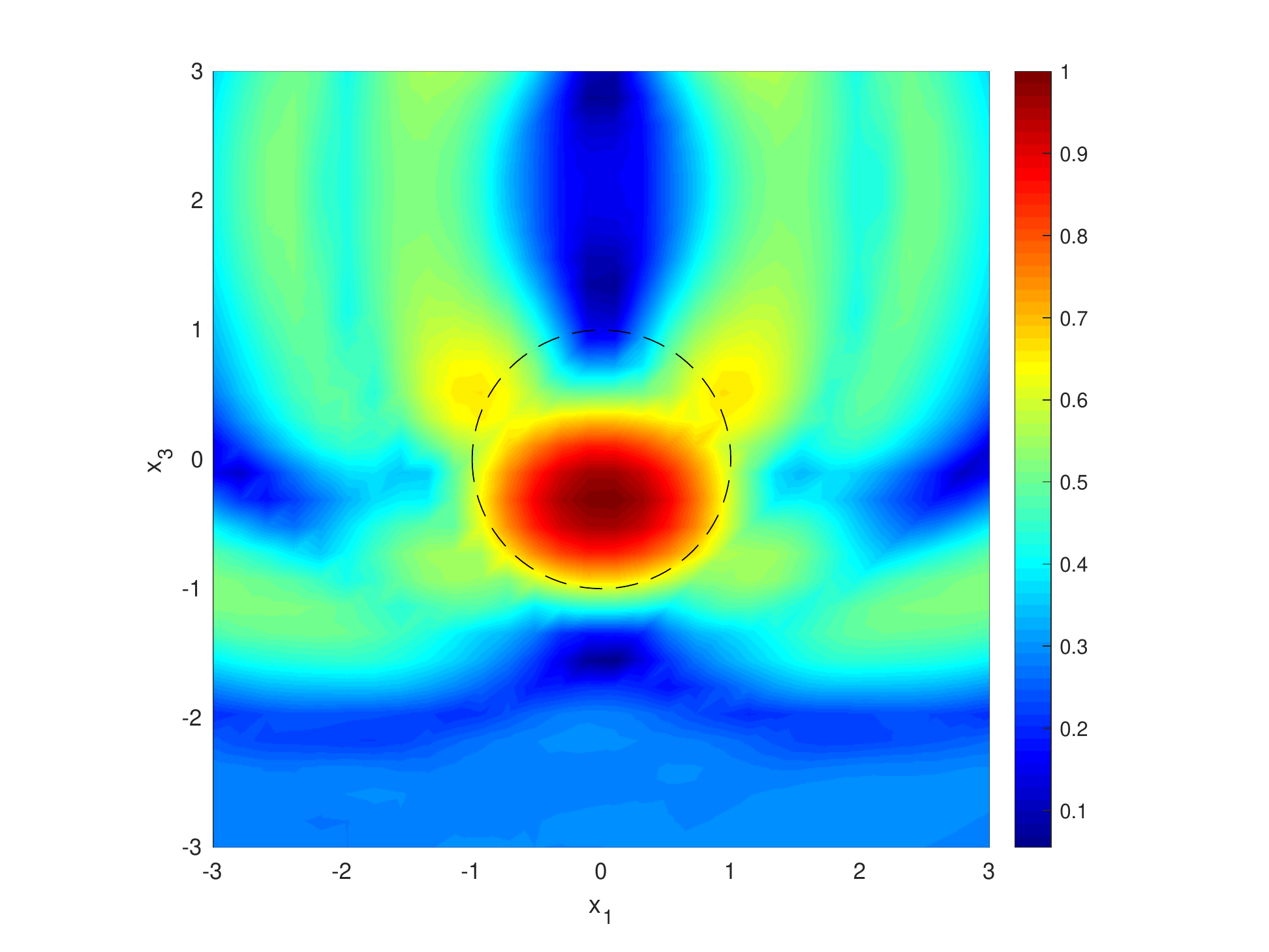}
\caption{{\bf Reconstructions of the ball  with $3$ measurement points located on the upper sphere.} Top left: exact ball. Top middle and top right: iso-surface view of the reconstruction with iso-values $8.5\times 10^{-1}$ and $8\times 10^{-1}$ respectively. Bottom left: $x_1x_2$ cross section view of the reconstruction. Bottom middle: $x_2x_3$ cross section view of the reconstruction. Bottom right: $x_1x_3$ cross section view of the reconstruction. }
\label{ball pt 3 upper}
\end{figure}

What may we reconstruct with a little bit more measurements points? Now we still consider the reconstruction of the same ball \eqref{numeric ball def} but with $3$ measurement points given by Table \ref{pt 3 table}, where these $3$ measurement points are on the upper half sphere. We plot both the three dimensional reconstruction and its cross-section views in Figure \ref{ball pt 3 upper}. It is observed that the location of ball is clearly indicated with only $3$ measurement points.
\begin{table}[htp]
\begin{center}
\begin{tabular}{||c||c|}
$\phi$ \qquad $\theta$\\
------------------------------------------------------------\\
  -180.0000000000000  \qquad 45.000000000000000     \\
  -90.00000000000000  \qquad  45.000000000000000     \\
  0.0000000000000000  \qquad 45.000000000000000     \\
\end{tabular}
\end{center}
\caption{$3$ measurement points: $(3 \sin \theta \cos \phi, 3 \sin \theta \cos \phi,3   \cos \theta)$. Angles are in degrees.}
\label{pt 3 table}
\end{table}%

 \begin{figure}[htbp]
  \centering
    \includegraphics[width=5cm]{figures/Ball_exact-eps-converted-to.pdf}
  \includegraphics[width=5cm]{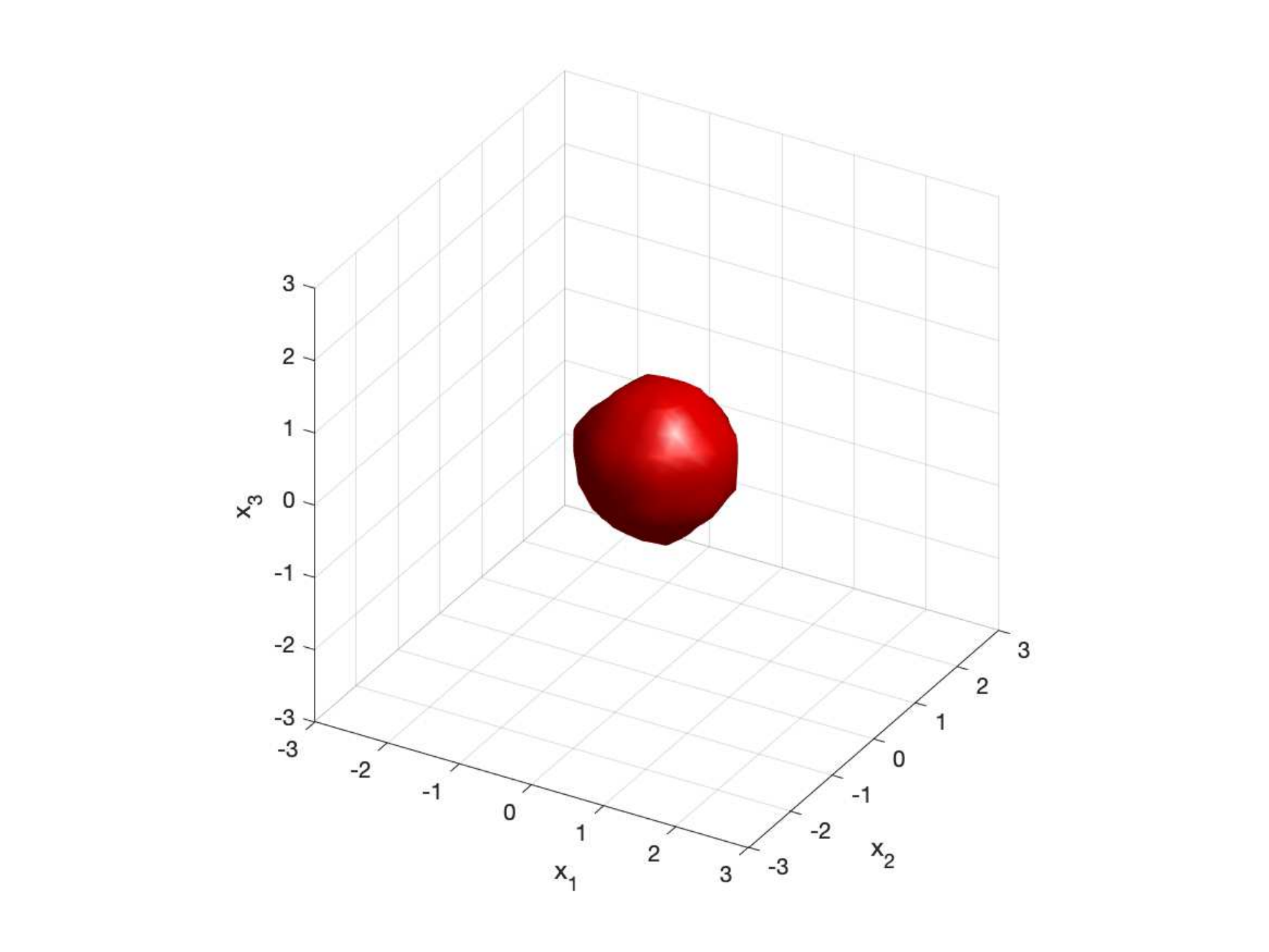}
    \includegraphics[width=5cm]{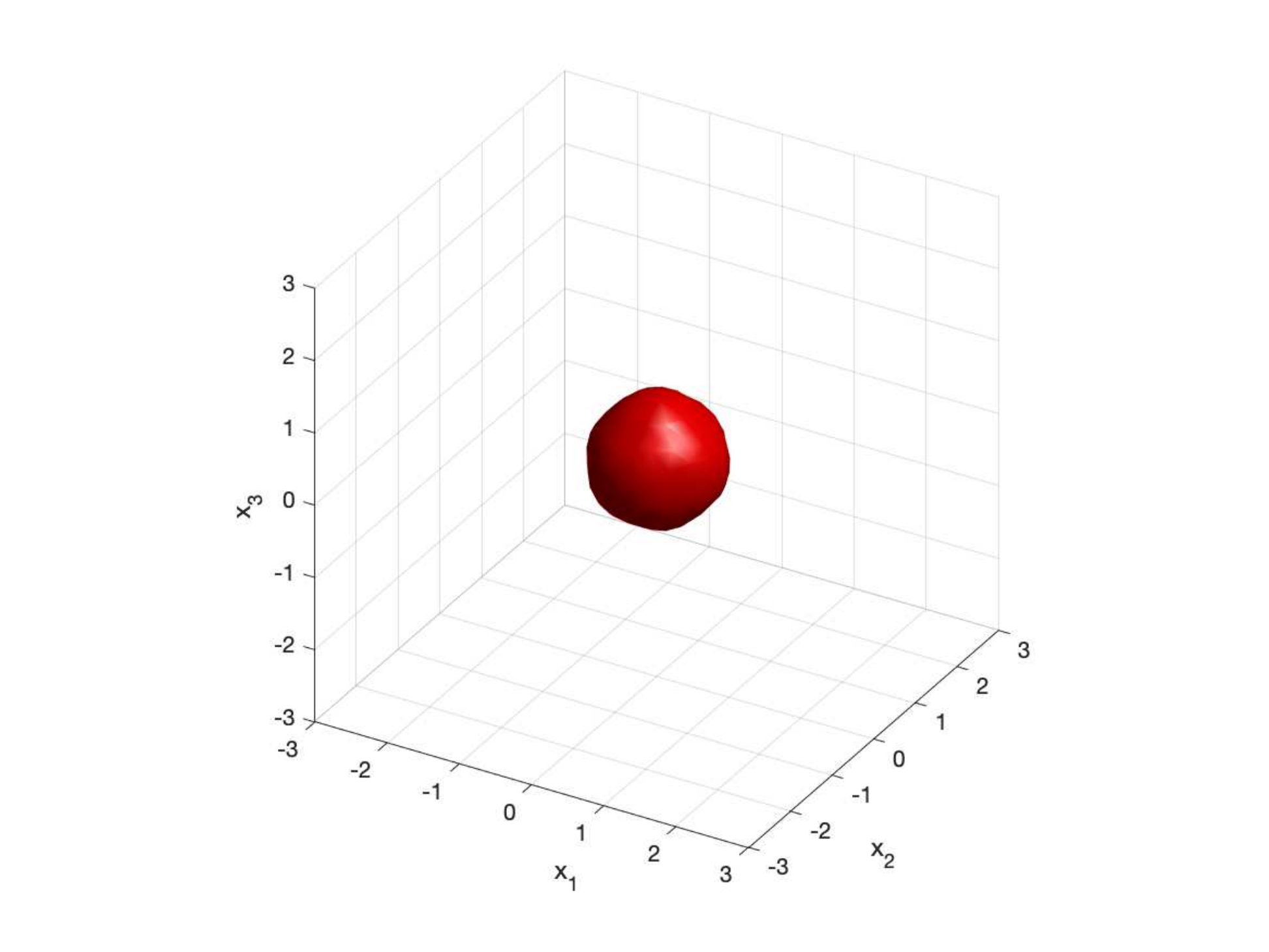}
    \includegraphics[width=5cm]{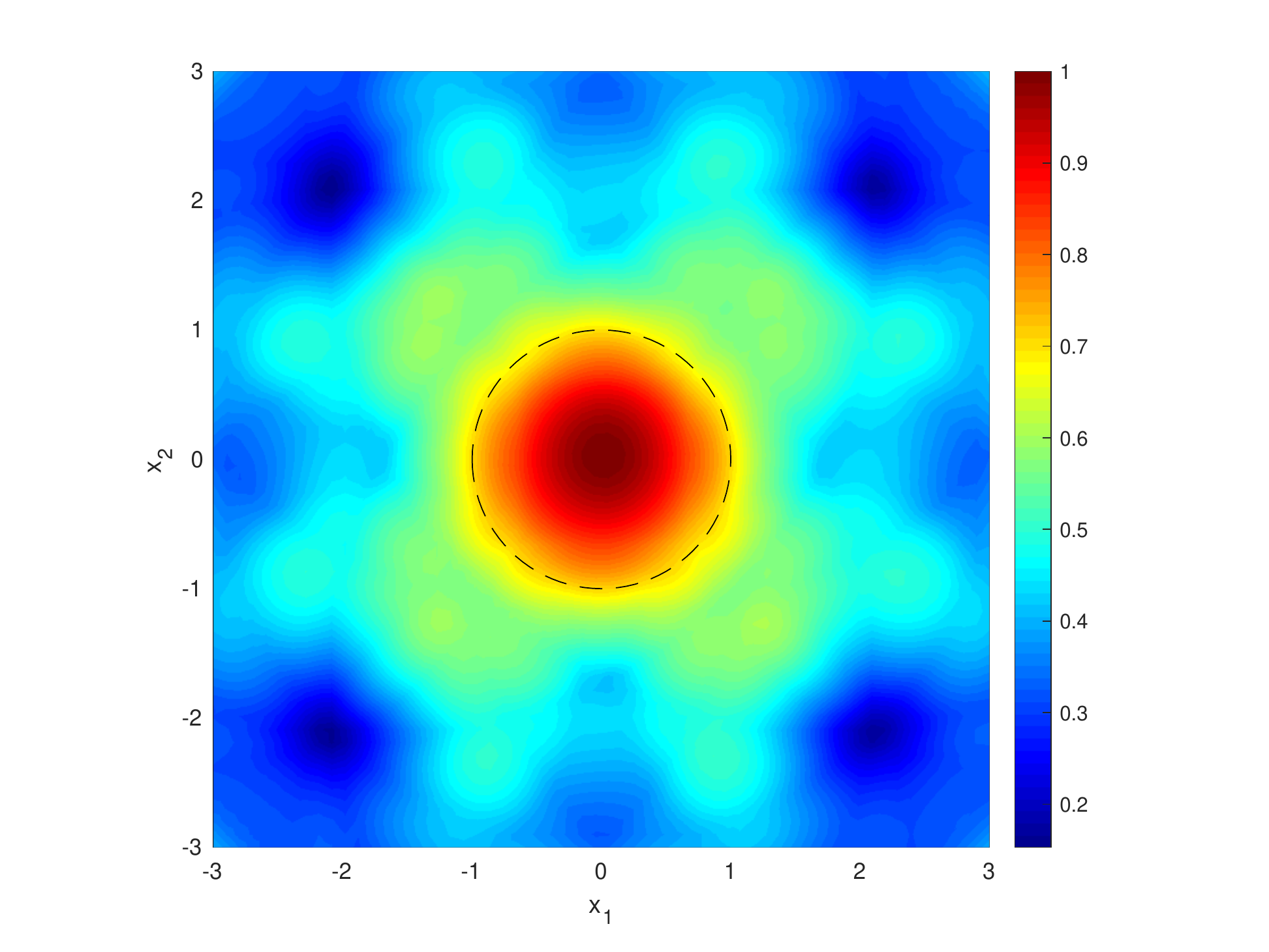}
     \includegraphics[width=5cm]{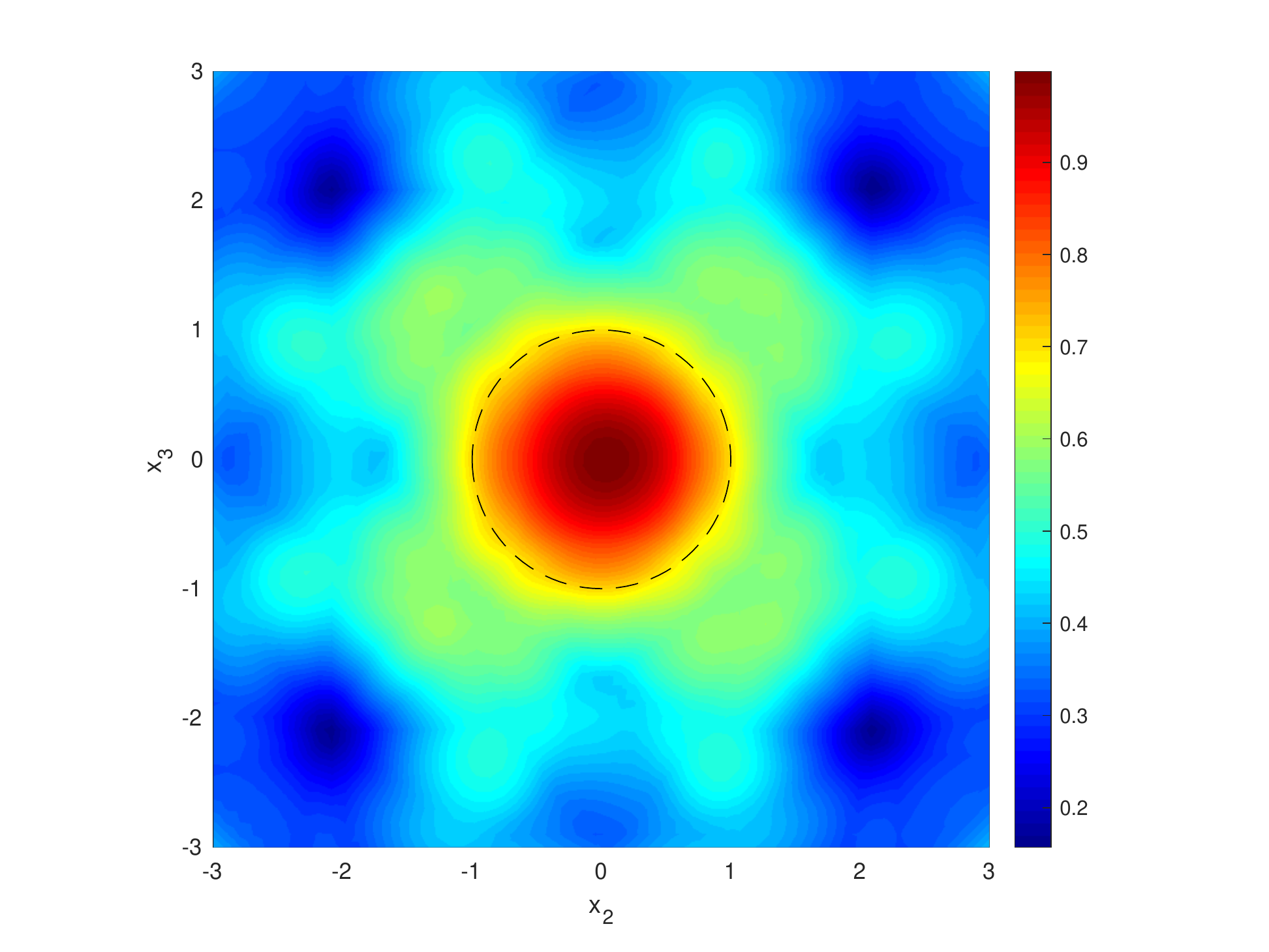}
          \includegraphics[width=5cm]{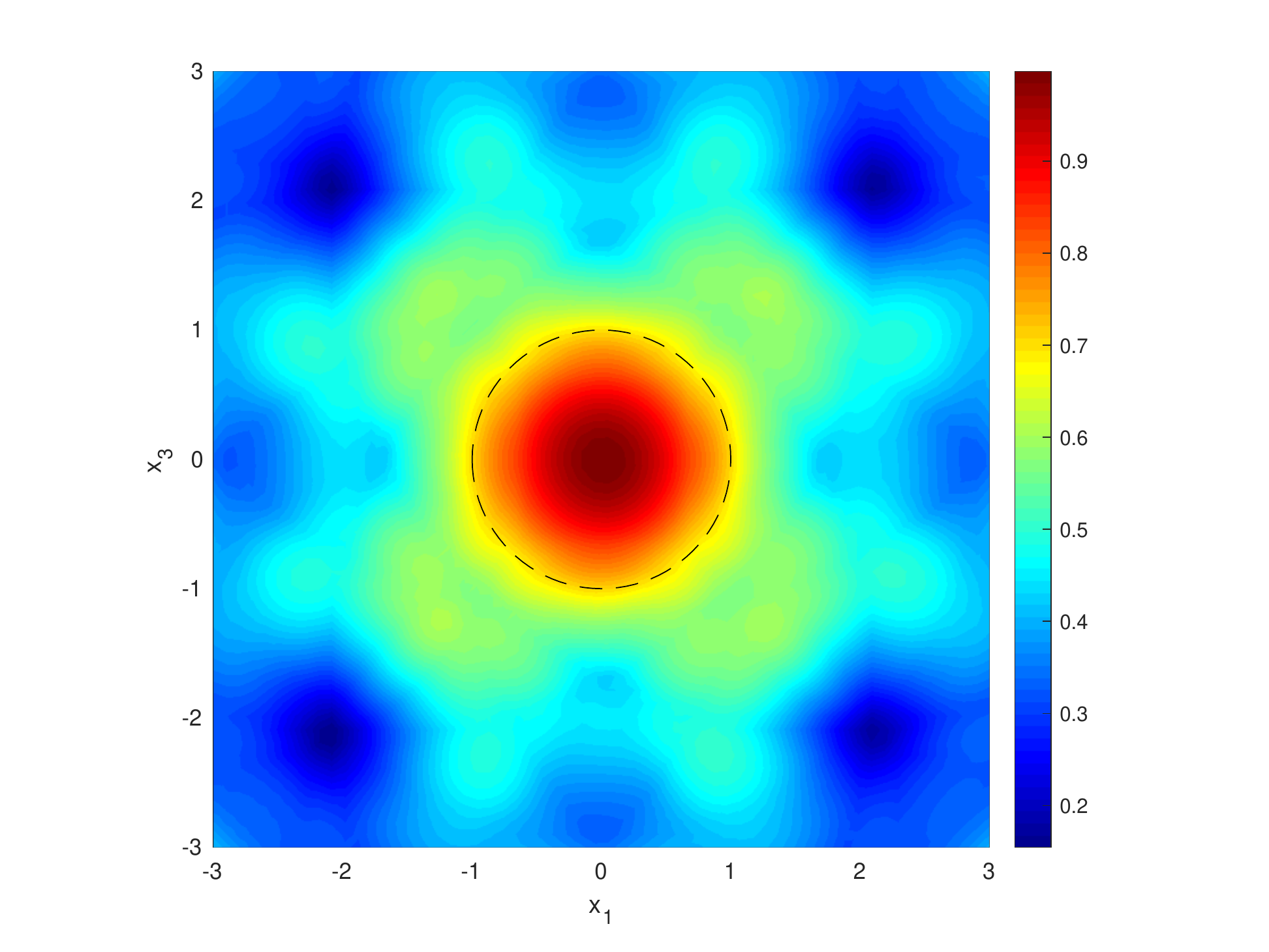}
\caption{{\bf Reconstructions of the ball  with $14$ measurement points.} Top left: exact ball. Top middle and top right: iso-surface view of the reconstruction with iso-values $7\times 10^{-1}$ and $7.5\times 10^{-1}$ respectively. Bottom left: $x_1x_2$ cross section view of the reconstruction. Bottom middle: $x_2x_3$ cross section view of the reconstruction. Bottom right: $x_1x_3$ cross section view of the reconstruction. }
\label{ball pt 14}
\end{figure}

To further shed light on the performance of our imaging algorithm, we consider the case when multiple (but still sparse) measurement points are available all around the source. In the following examples, we consider $14$ measurement points which are roughly equally distributed on the measurement sphere. We remark that only $14\times11=154$ measurements (in three dimensional space and frequency) are used in this case. The polar coordinates of the $14$ measurement points are given in Table \ref{pt 14 table}.

Figure \ref{ball pt 14} shows the reconstruction of the ball \eqref{numeric ball def} using multiple but sparse measurements. Compared with the reconstruction using $3$ measurement points, it is observed that both the location and shape of the ball are better reconstructed.

\begin{table}[htp]
\begin{center}
\begin{tabular}{||c||c|}
$\phi$ \qquad $\theta$\\
------------------------------------------------------------\\
   0.0000000000000000   \qquad  90.0000000000000000     \\
 180.00000000000000  \qquad  90.0000000000000000     \\
  90.000000000000000   \qquad 90.0000000000000000     \\
 -90.00000000000000  \qquad  90.0000000000000000     \\
  90.000000000000000  \qquad   0.00000000000000000     \\
  90.000000000000000  \qquad 180.000000000000000     \\
  45.000000000000000  \qquad  54.7356103172453460     \\
  45.000000000000000  \qquad 125.264389682754654     \\
 -45.00000000000000  \qquad  54.7356103172453460     \\
 -45.00000000000000  \qquad 125.264389682754654     \\
 135.00000000000000 \qquad   54.7356103172453460     \\
 135.00000000000000 \qquad  125.264389682754654     \\
-135.0000000000000  \qquad  54.7356103172453460     \\
-135.0000000000000  \qquad 125.264389682754654     \\
\end{tabular}
\end{center}
\caption{$14$ measurement points: $(r_m \sin \theta \cos \phi, r_m \sin \theta \cos \phi,r_m   \cos \theta)$. Angles are in degrees.}
\label{pt 14 table}
\end{table}%

\begin{figure}[htbp]
  \centering
  \includegraphics[width=5cm]{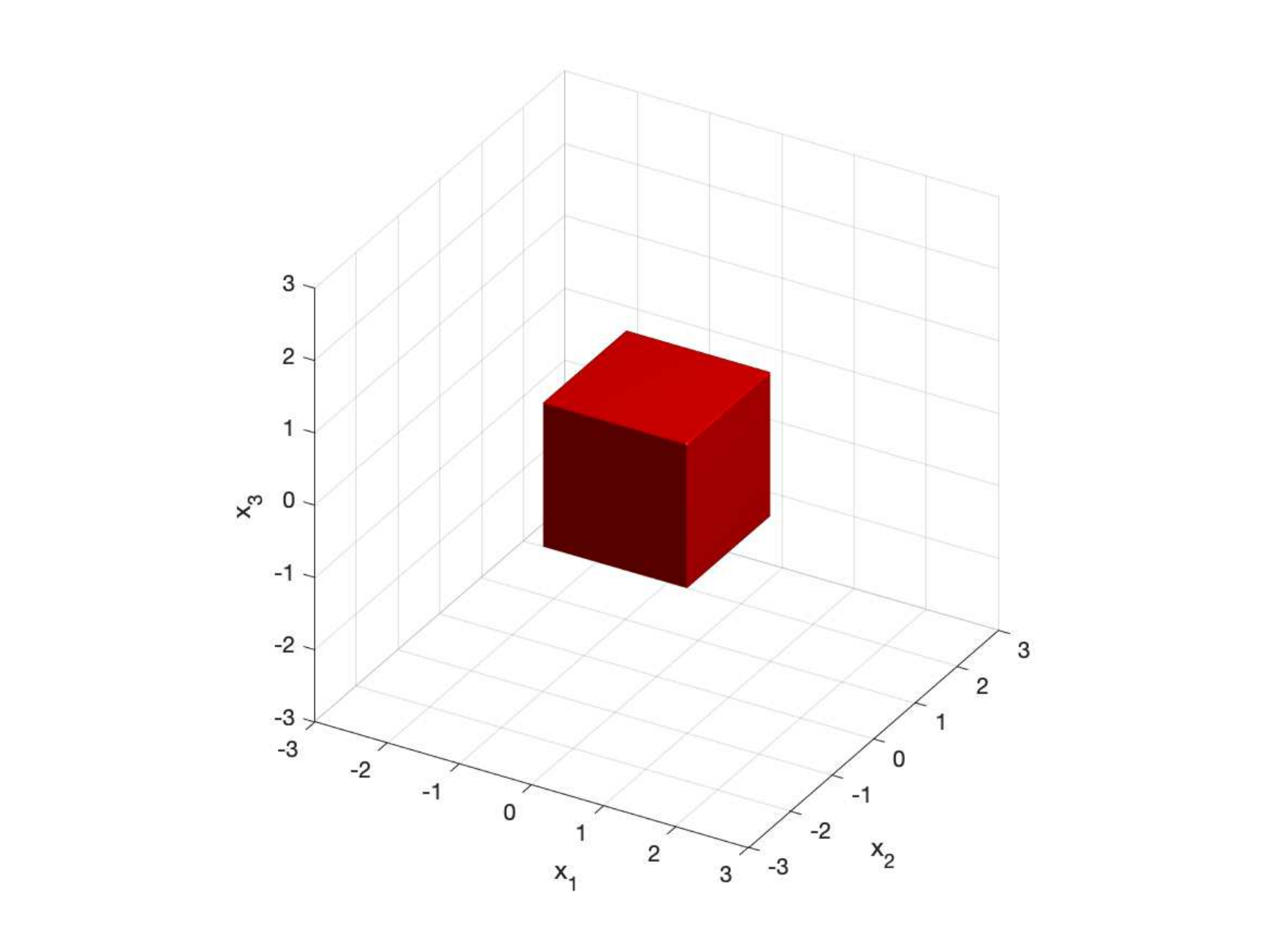}
      \includegraphics[width=5cm]{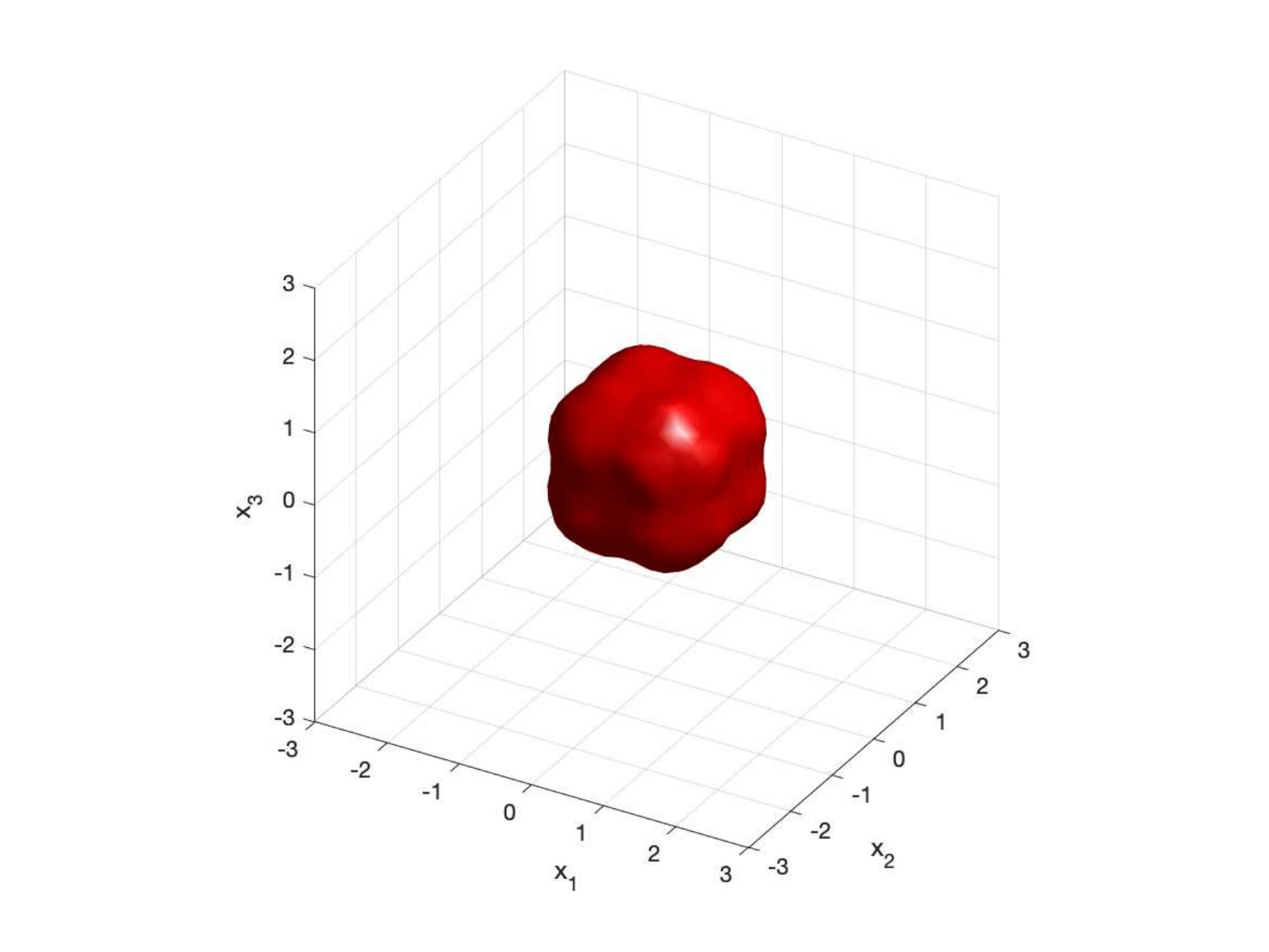}
      \includegraphics[width=5cm]{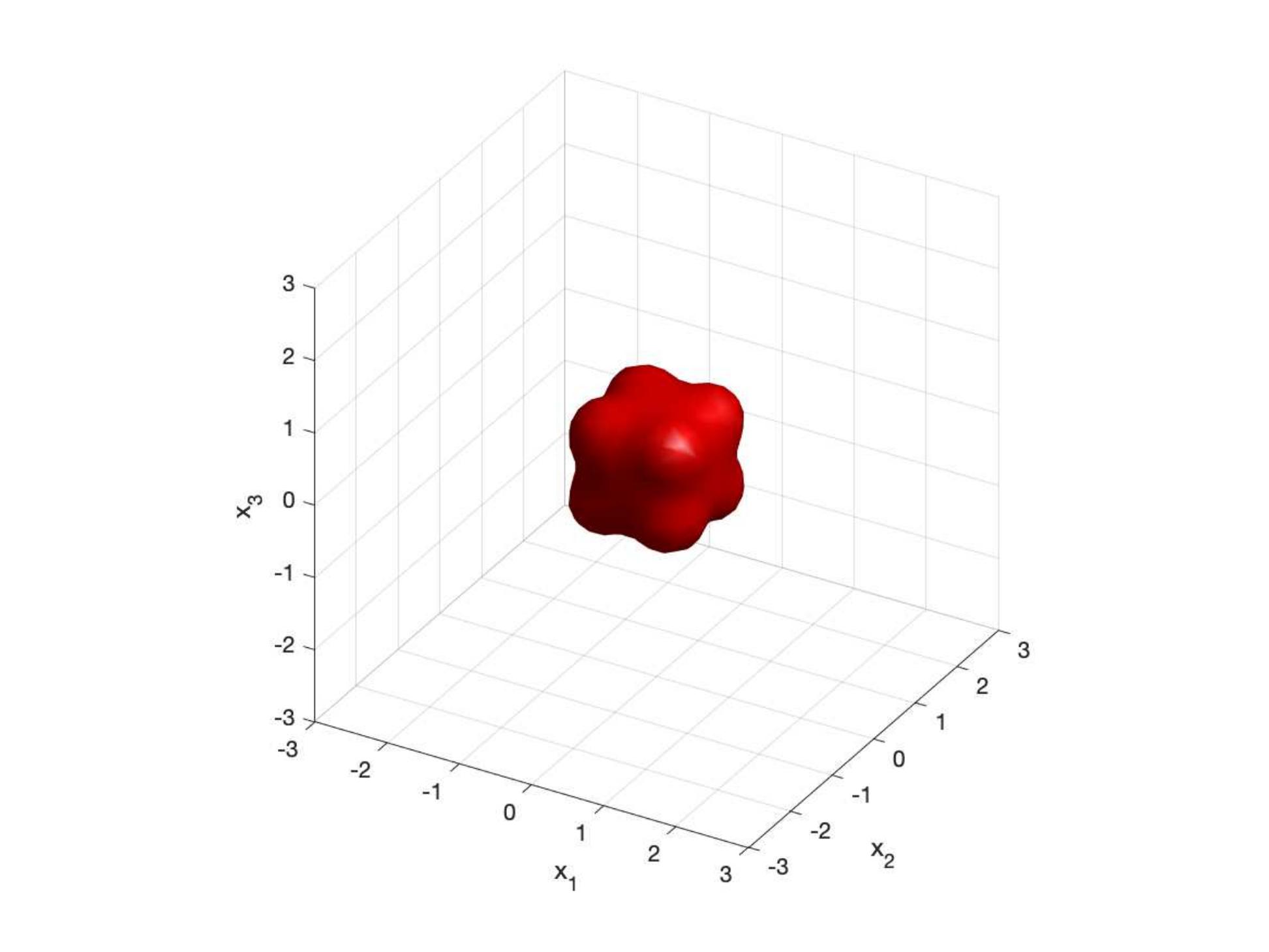}
    \includegraphics[width=5cm]{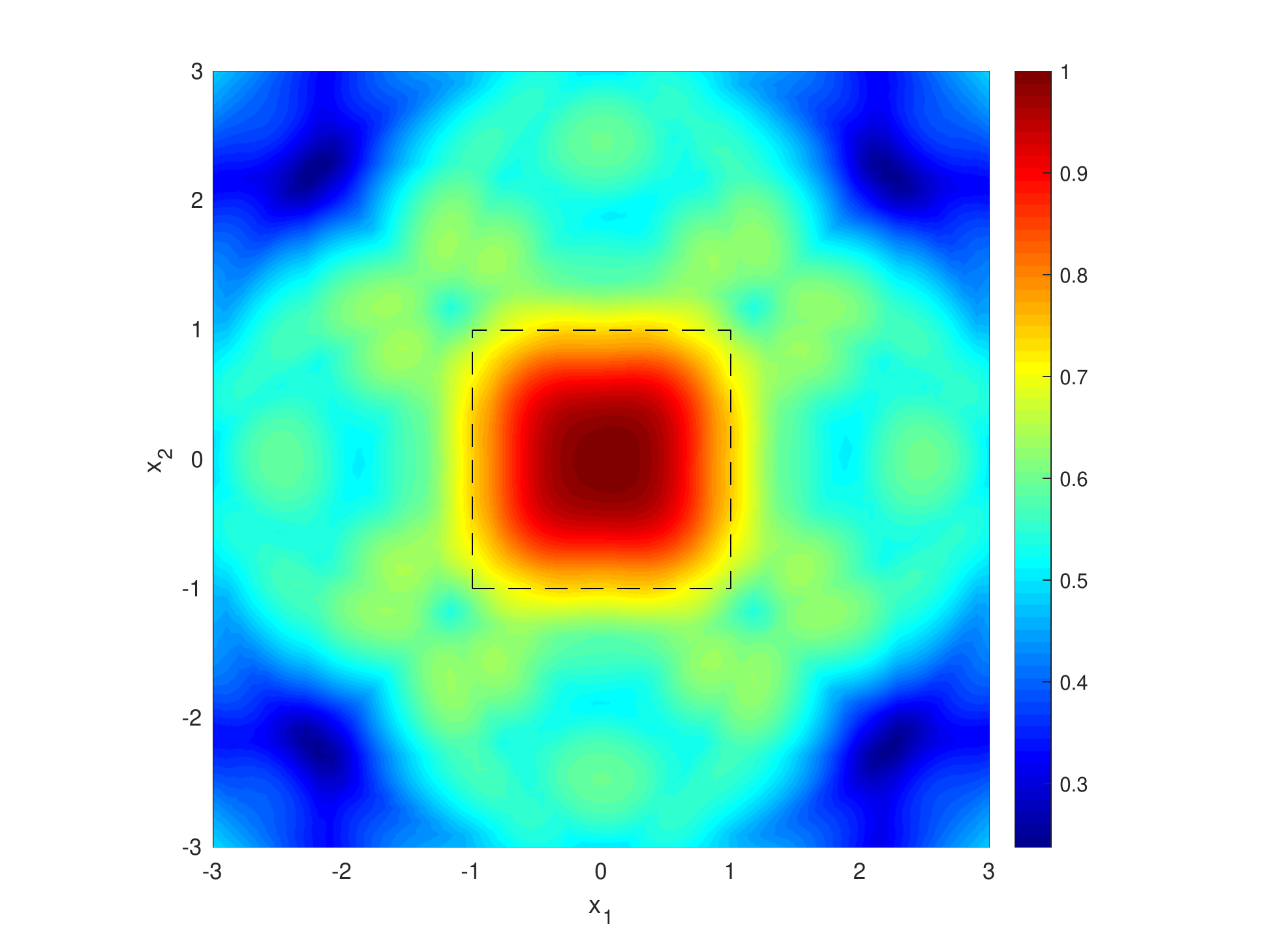}
     \includegraphics[width=5cm]{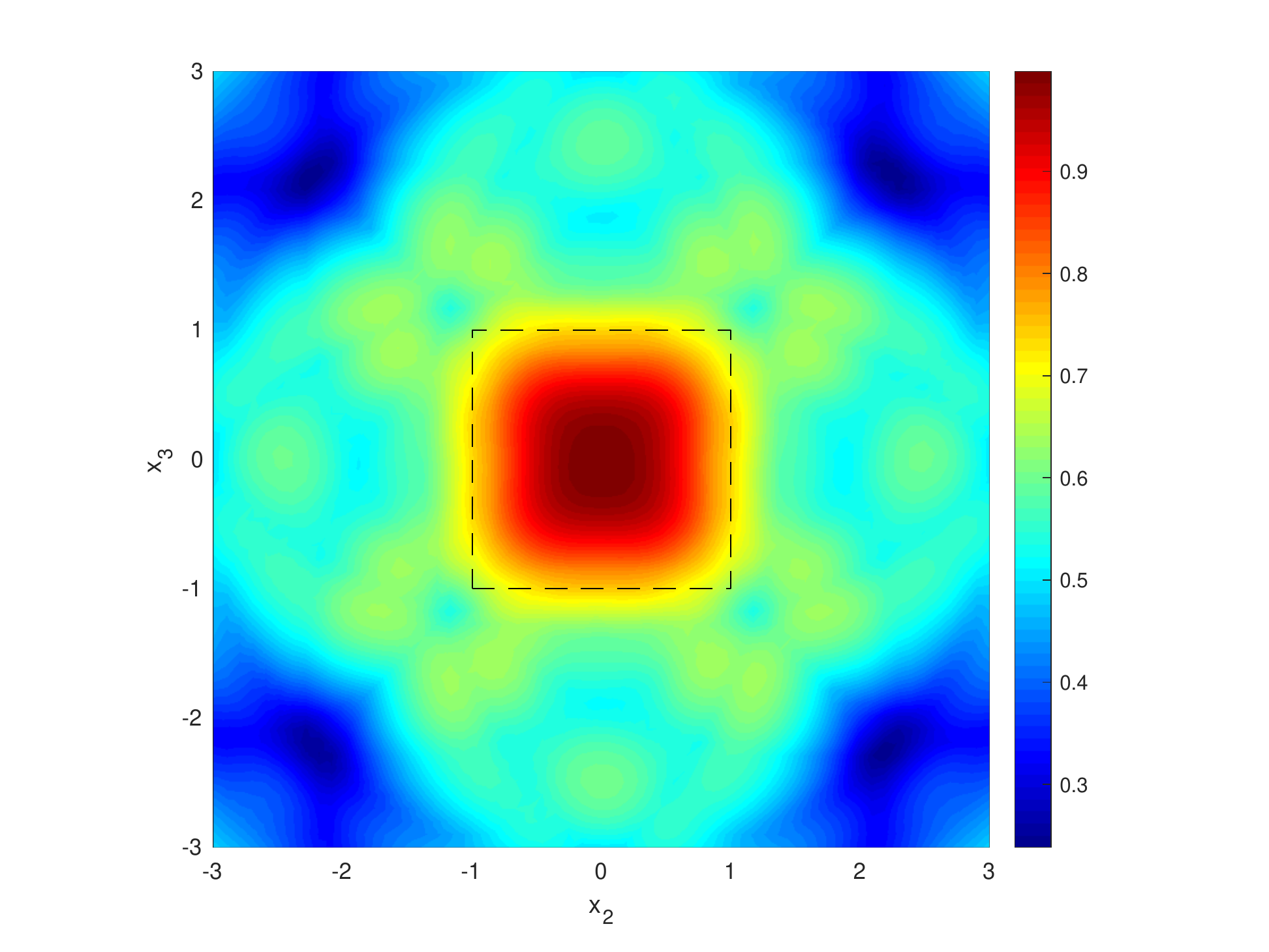}
          \includegraphics[width=5cm]{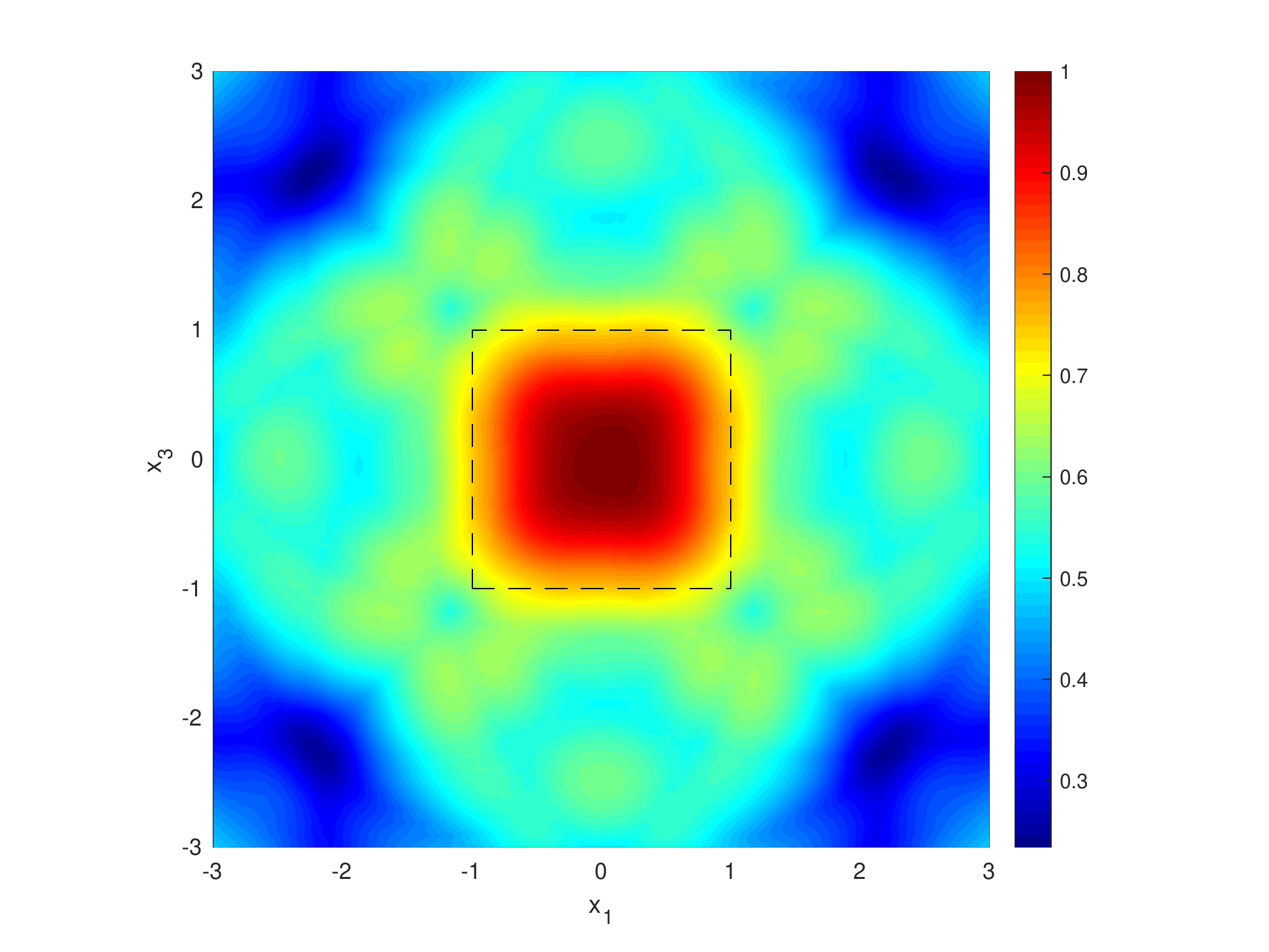}
\caption{{\bf Reconstructions of the cube with $14$ measurement points.} Top left: exact cube. Top middle and top right: iso-surface view of the reconstruction with iso-values $7\times 10^{-1}$ and $8\times 10^{-1}$ respectively. Bottom left: $x_1x_2$ cross section view of the reconstruction. Bottom middle: $x_2x_3$ cross section view of the reconstruction. Bottom right: $x_1x_3$ cross section view of the reconstruction. }
\label{cube pt 14}
\end{figure}

\begin{figure}[htbp]
  \centering
  \includegraphics[width=5cm]{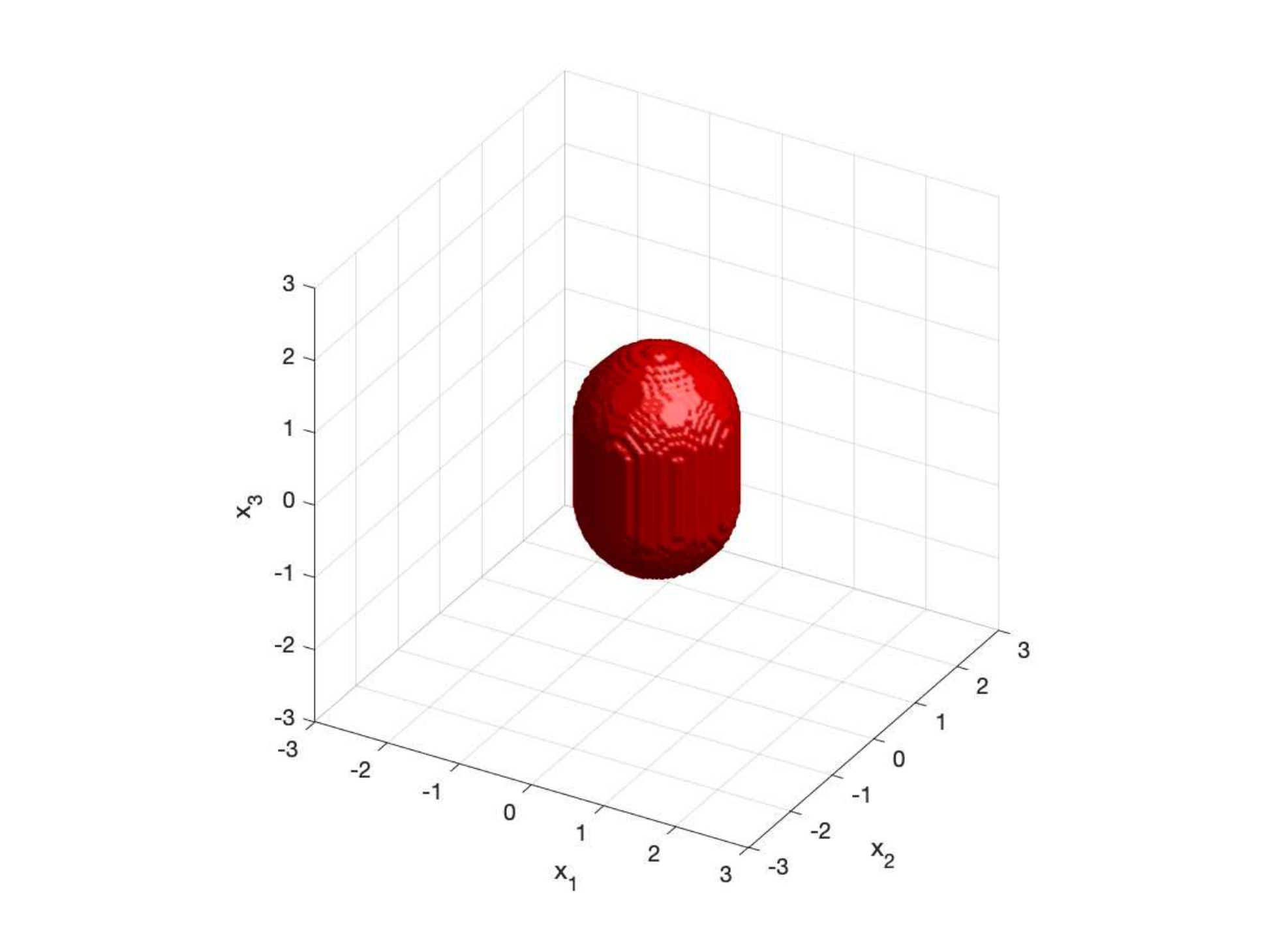}
    \includegraphics[width=5cm]{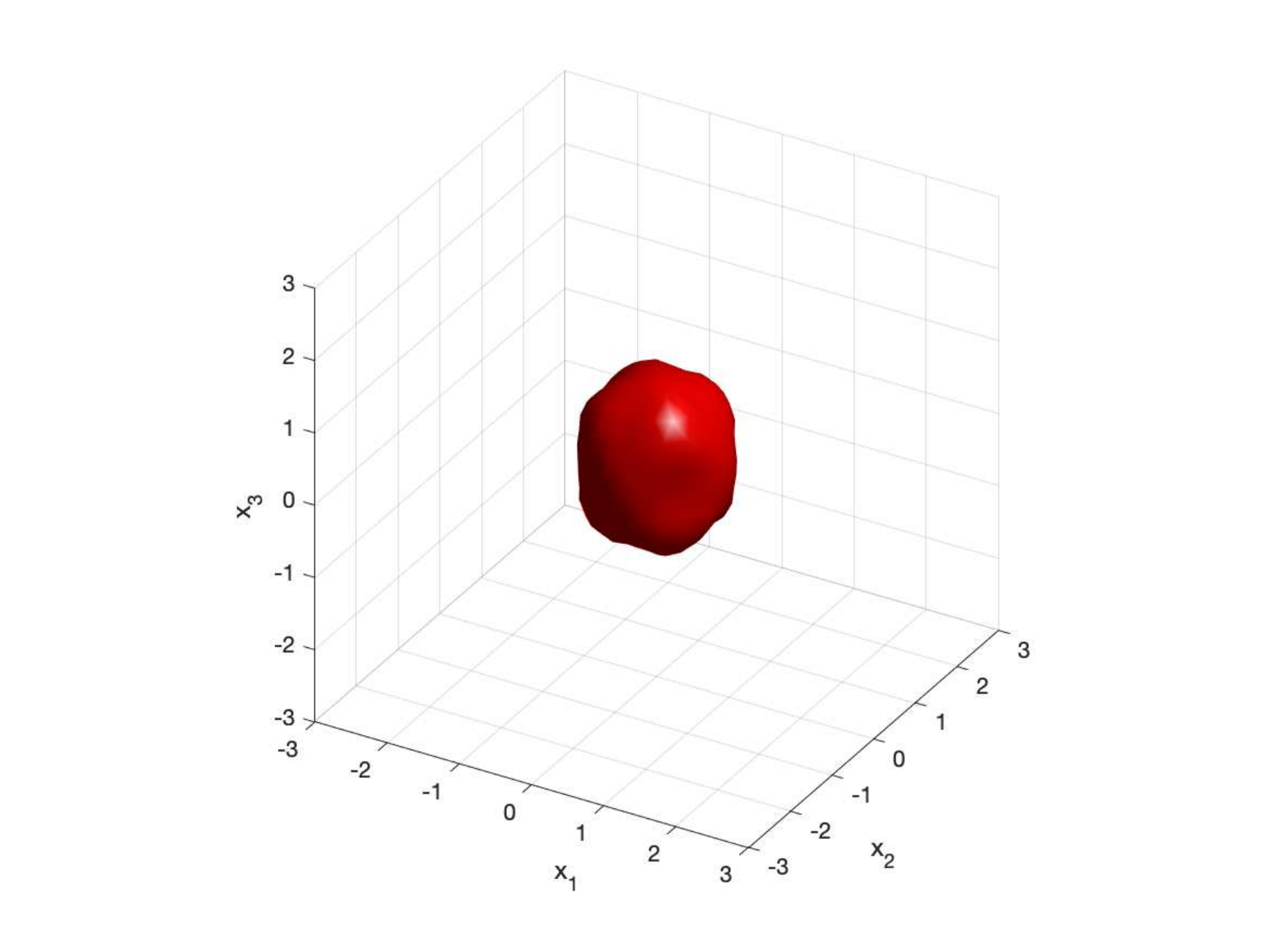}
    \includegraphics[width=5cm]{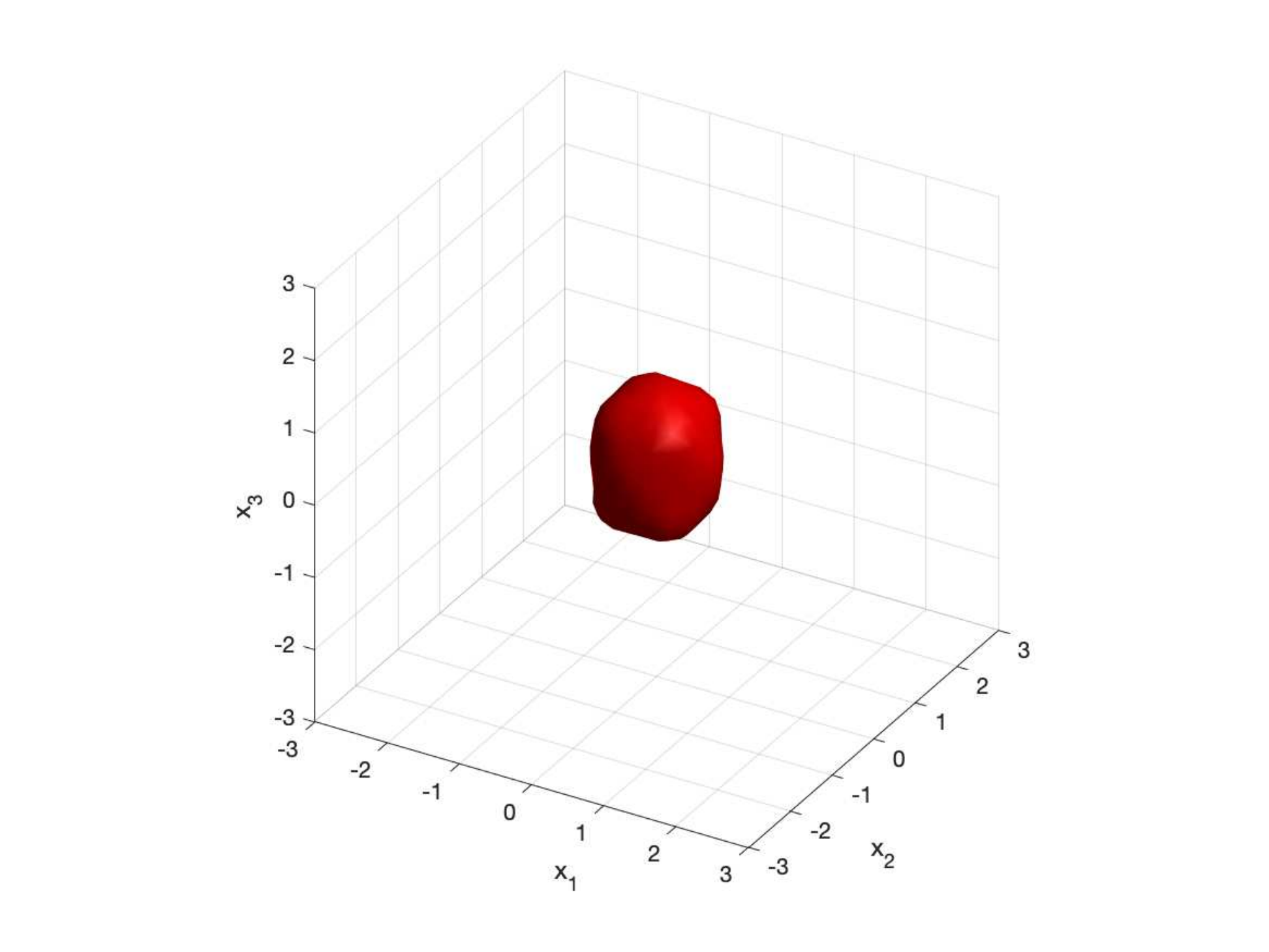}
        \includegraphics[width=5cm]{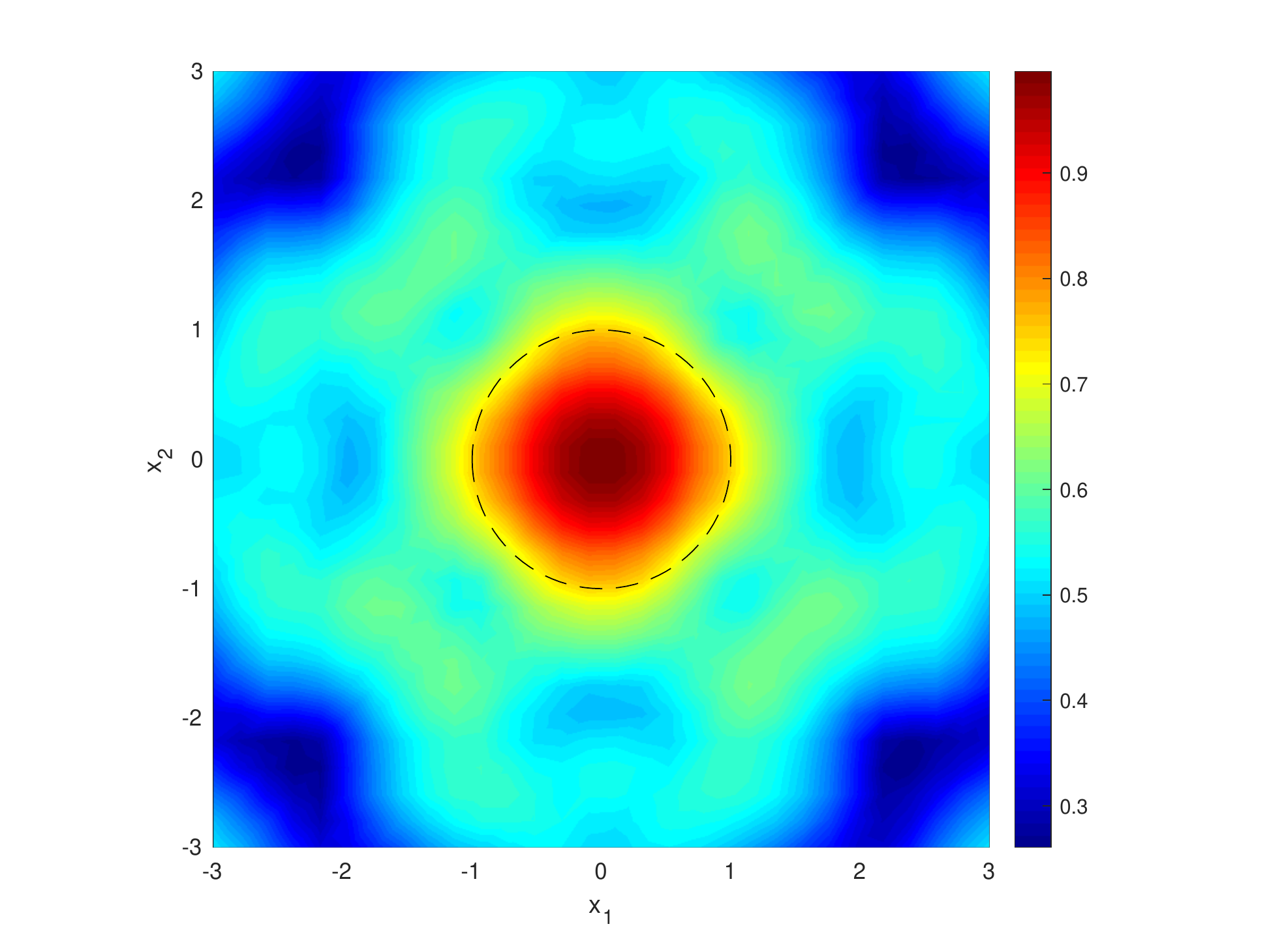}
        \includegraphics[width=5cm]{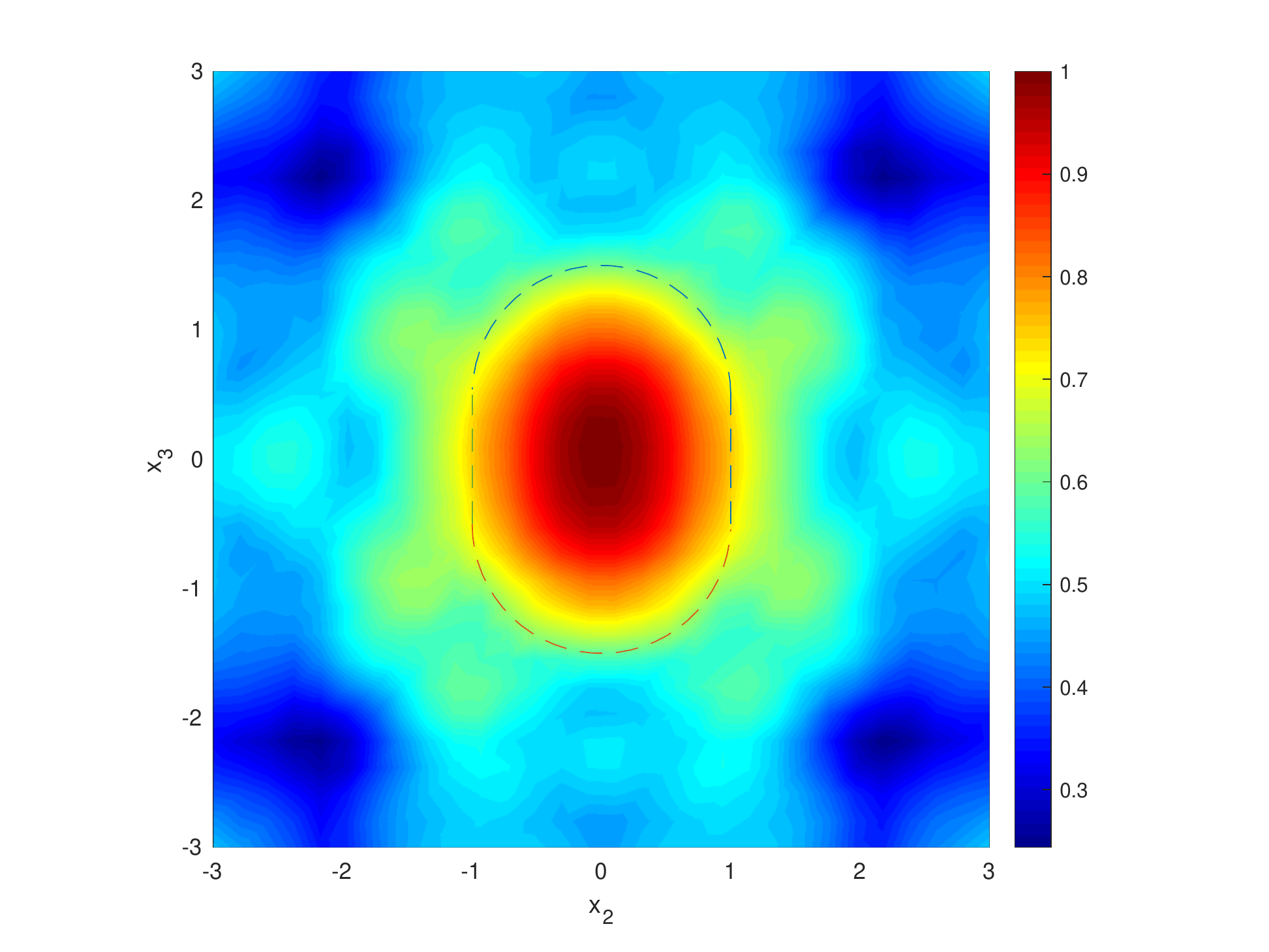}
        \includegraphics[width=5cm]{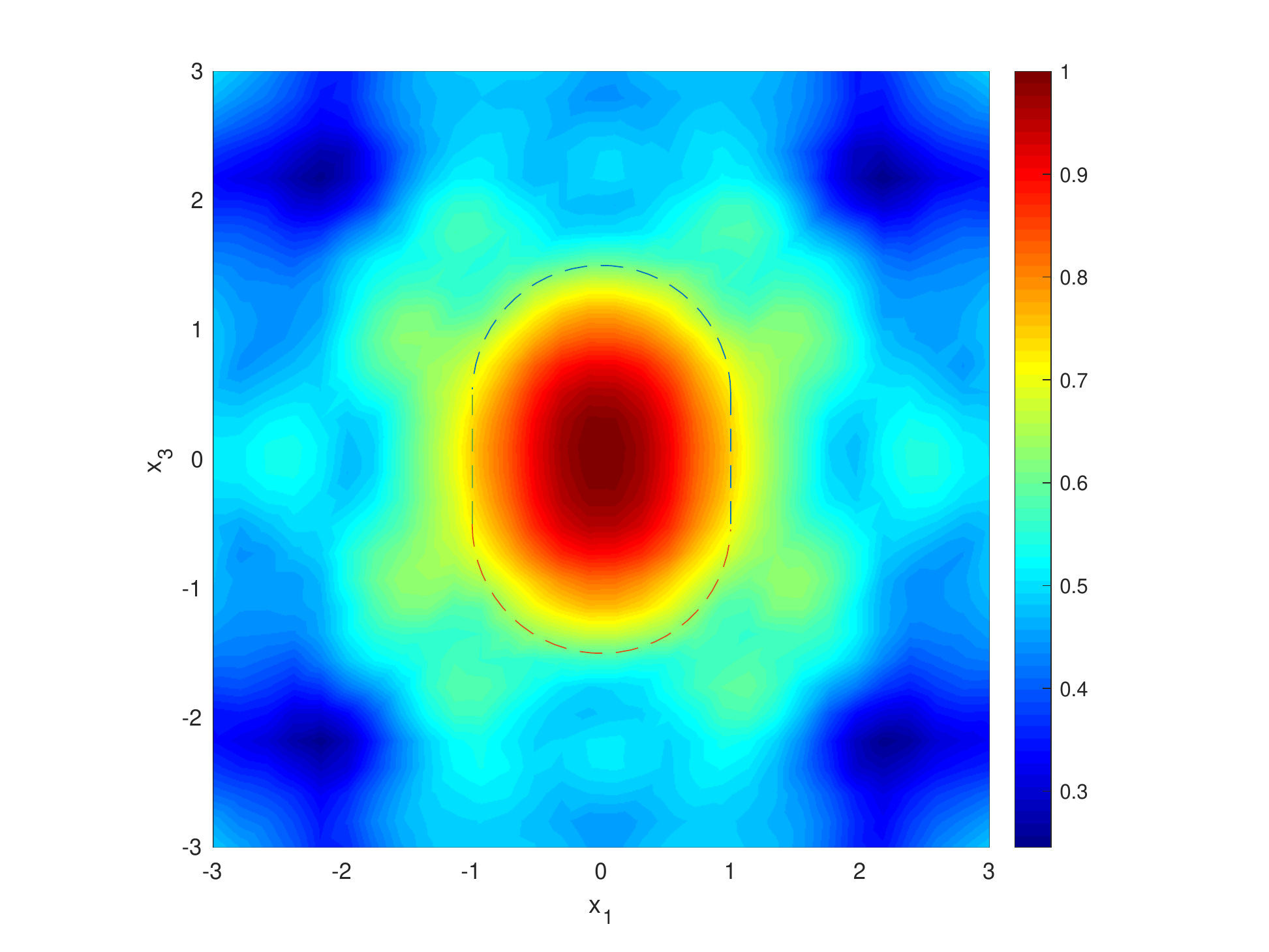}
\caption{{\bf Reconstructions of the rounded cylinder with $14$ measurement points.} Top left: exact rounded cylinder . Top middle and top right: iso-surface view of the reconstruction with iso-values $7.5\times 10^{-1}$ and $8\times 10^{-1}$ respectively. Bottom left: $x_1x_2$ cross section view of the reconstruction. Bottom middle: $x_2x_3$ cross section view of the reconstruction. Bottom right: $x_1x_3$ cross section view of the reconstruction.  }
\label{cr pt 14}
\end{figure}

\begin{figure}[htbp]
  \centering
  \includegraphics[width=5cm]{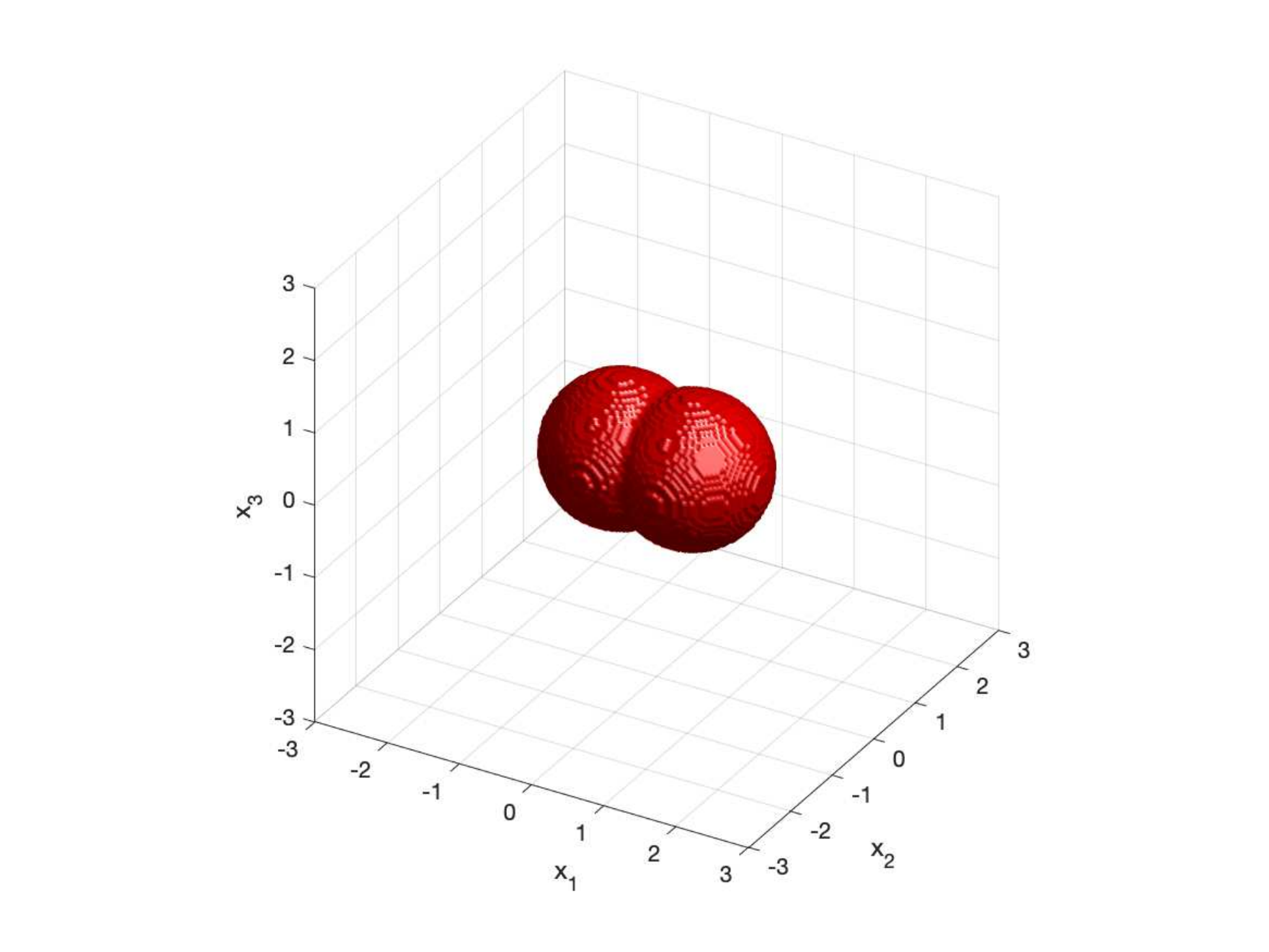}
    \includegraphics[width=5cm]{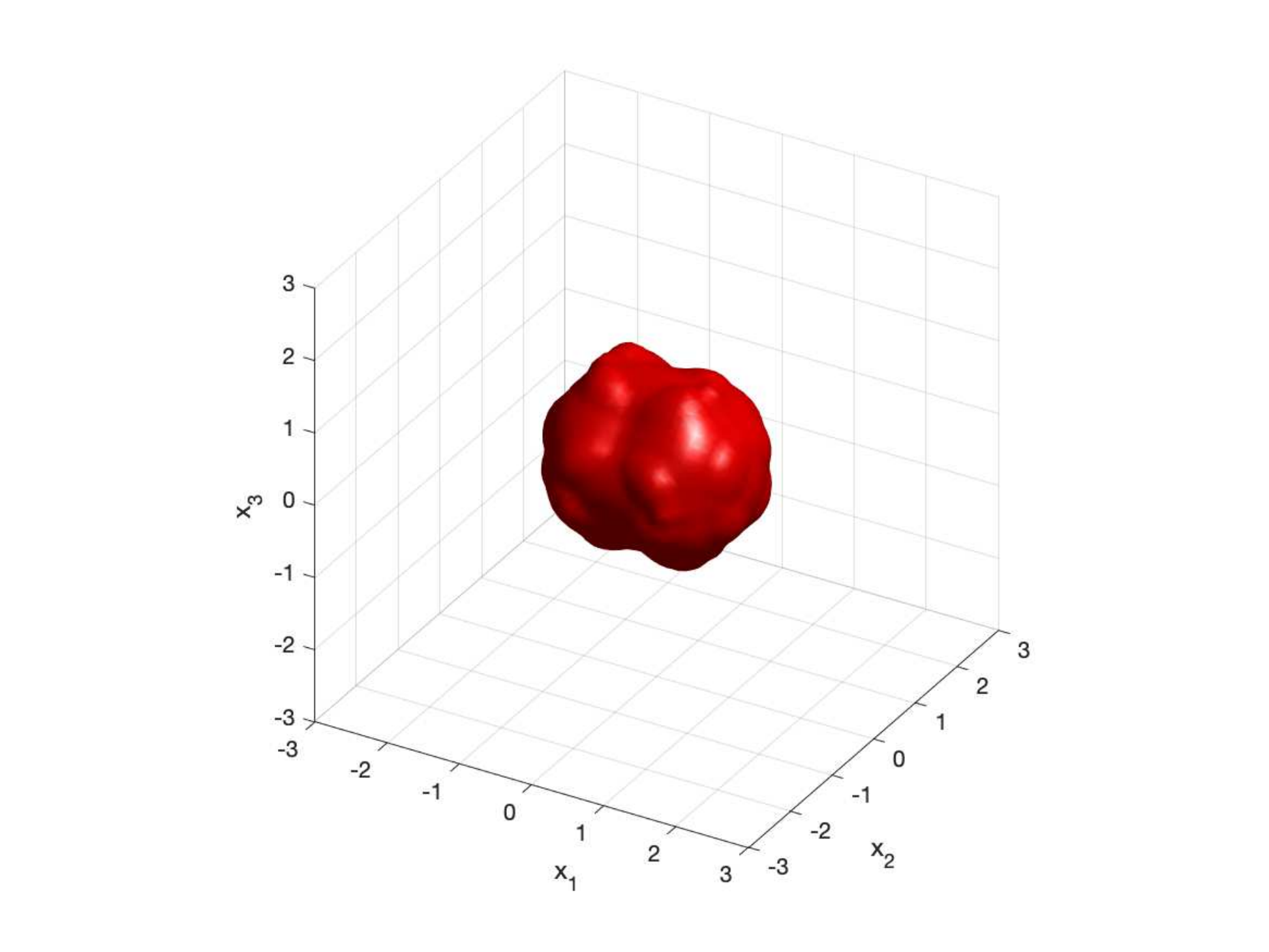}
    \includegraphics[width=5cm]{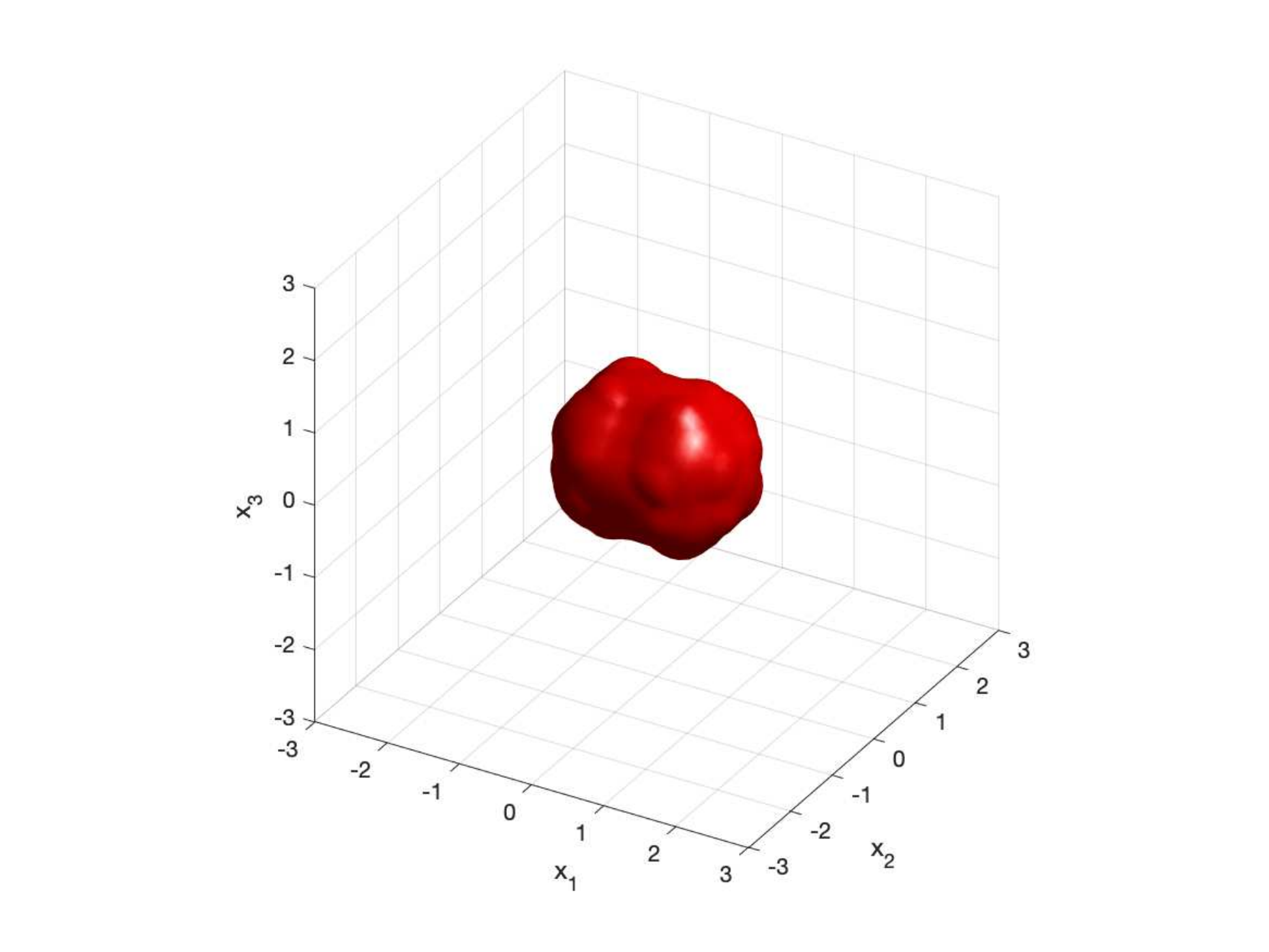}
        \includegraphics[width=5cm]{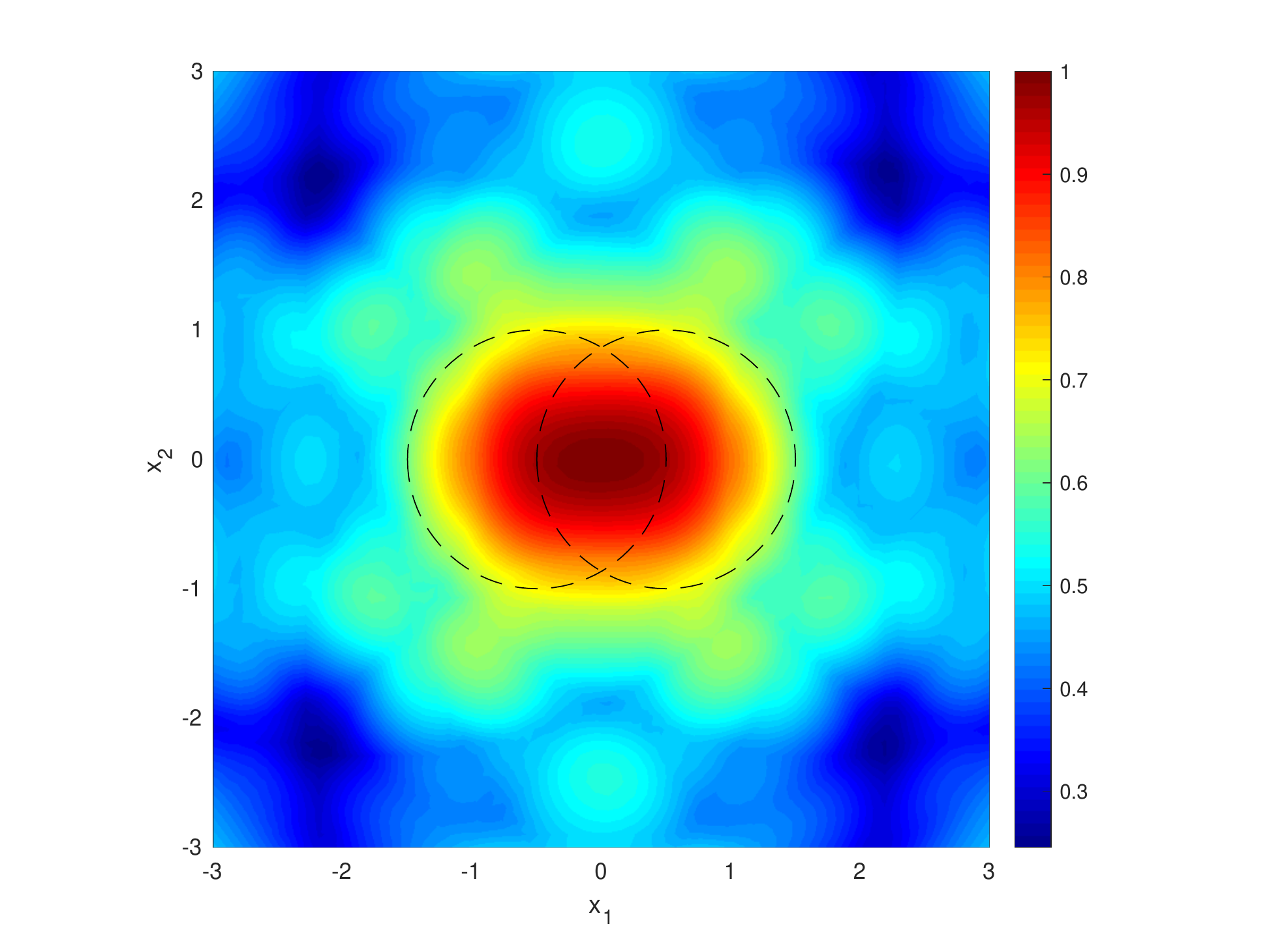}
        \includegraphics[width=5cm]{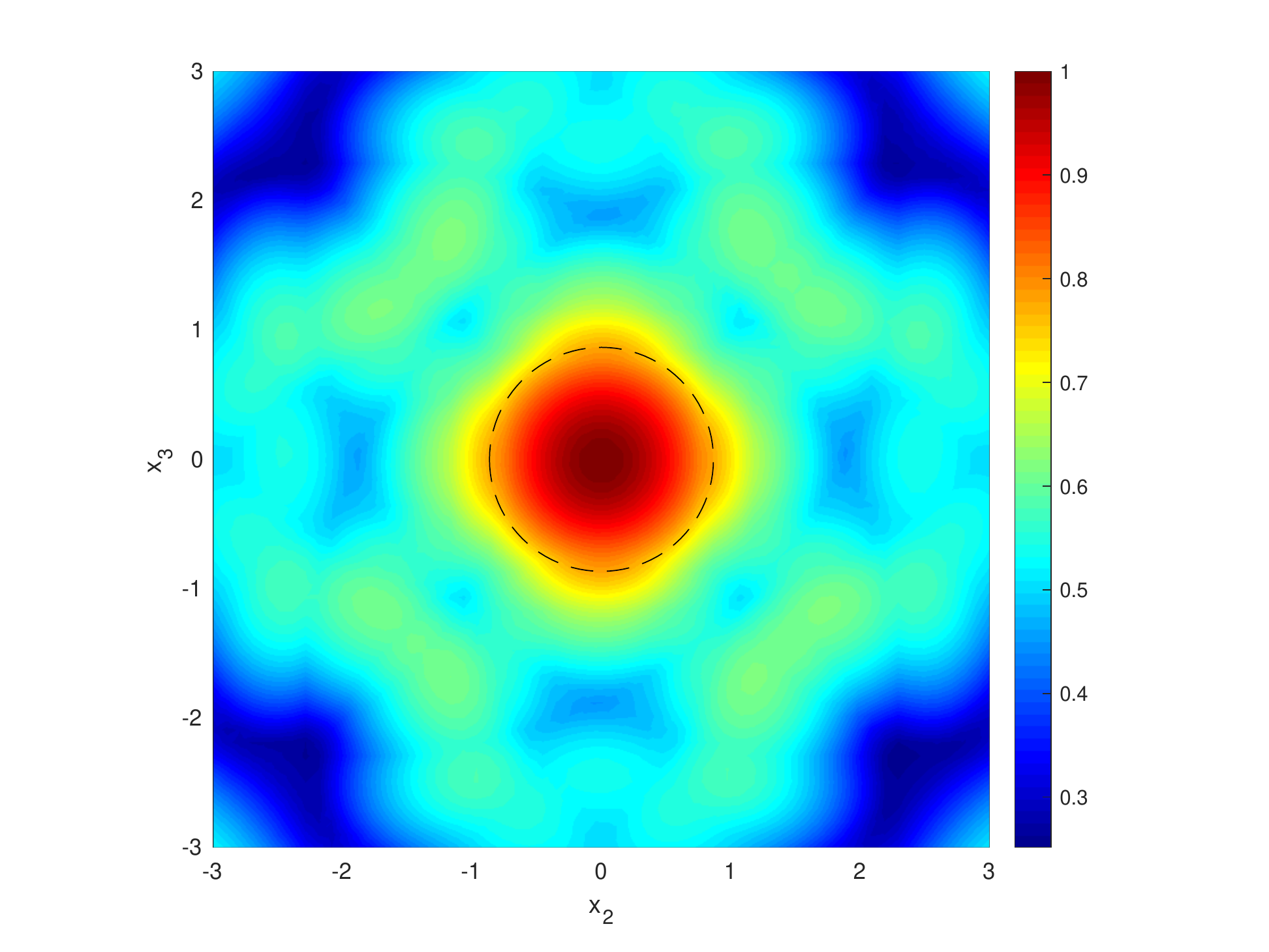}
        \includegraphics[width=5cm]{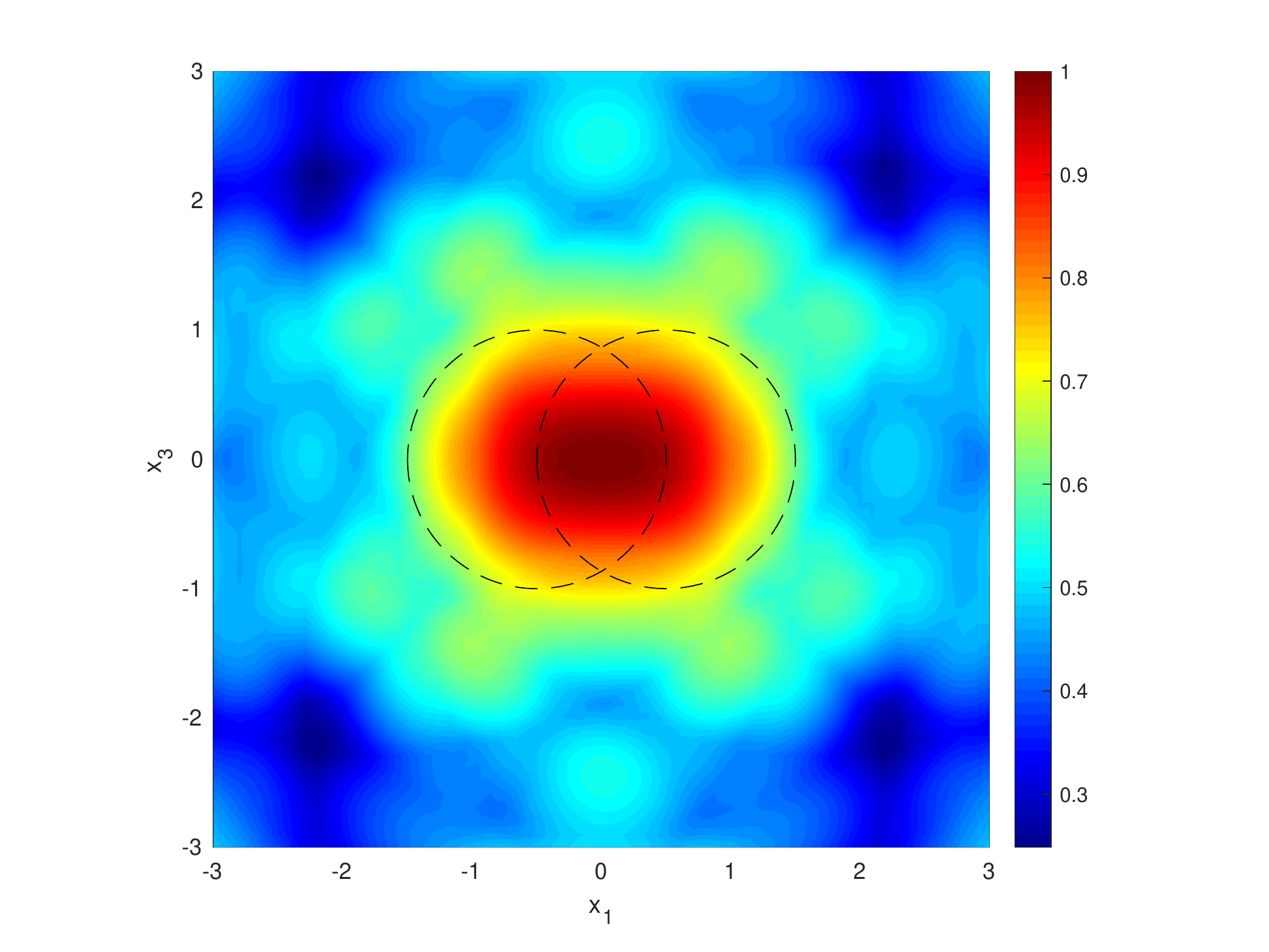}
\caption{{\bf Reconstructions of the peanut with $14$ measurement points.} Top left: exact peanut. Top middle and top right: iso-surface view of the reconstruction with iso-values $7.5\times 10^{-1}$ and $8\times 10^{-1}$ respectively. Bottom left: $x_1x_2$ cross section view of the reconstruction. Bottom middle: $x_2x_3$ cross section view of the reconstruction. Bottom right: $x_1x_3$ cross section view of the reconstruction.   }
\label{peanut pt 14}
\end{figure}

\begin{figure}[htbp]
  \centering
  \includegraphics[width=5cm]{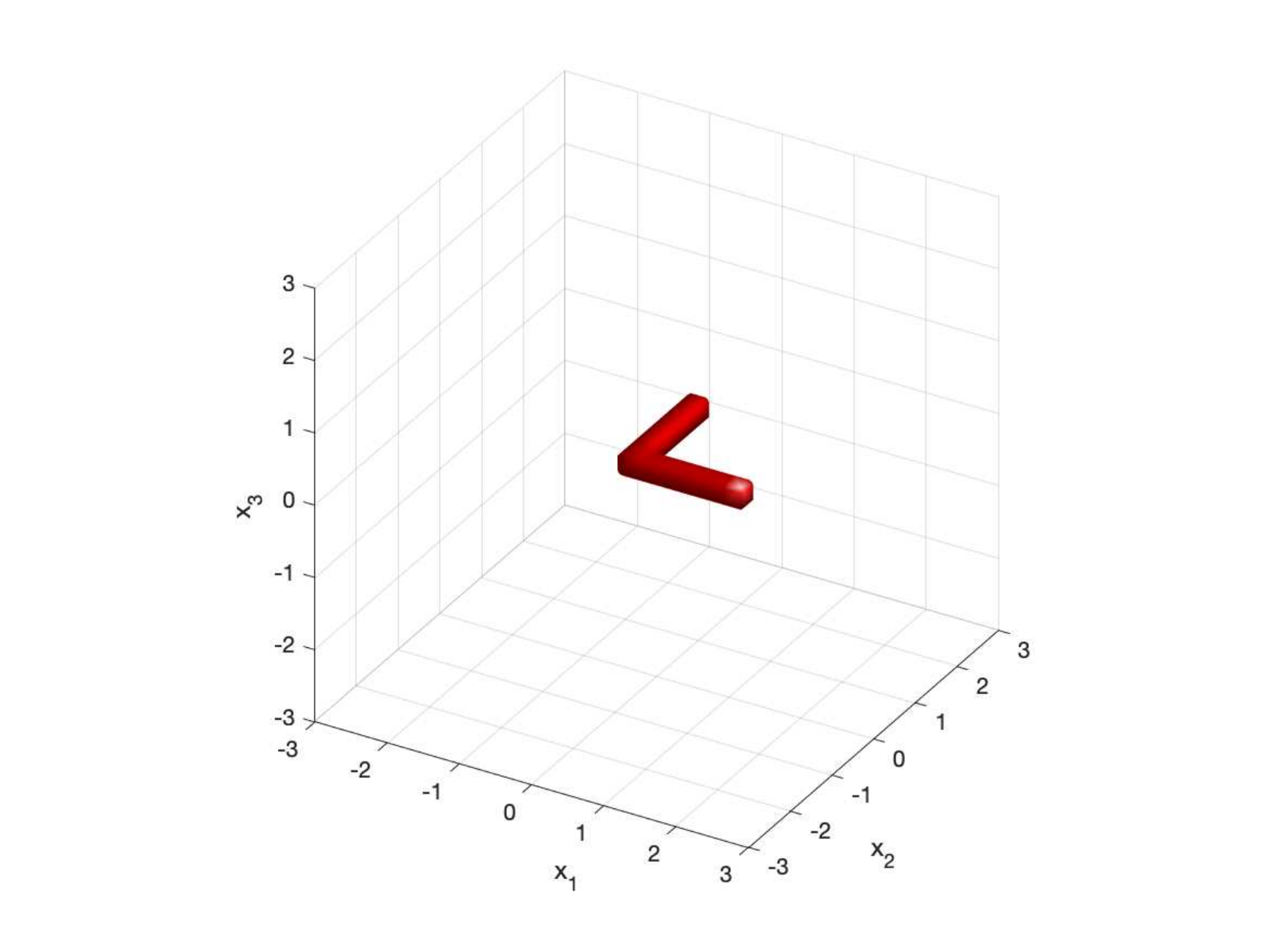}
    \includegraphics[width=5cm]{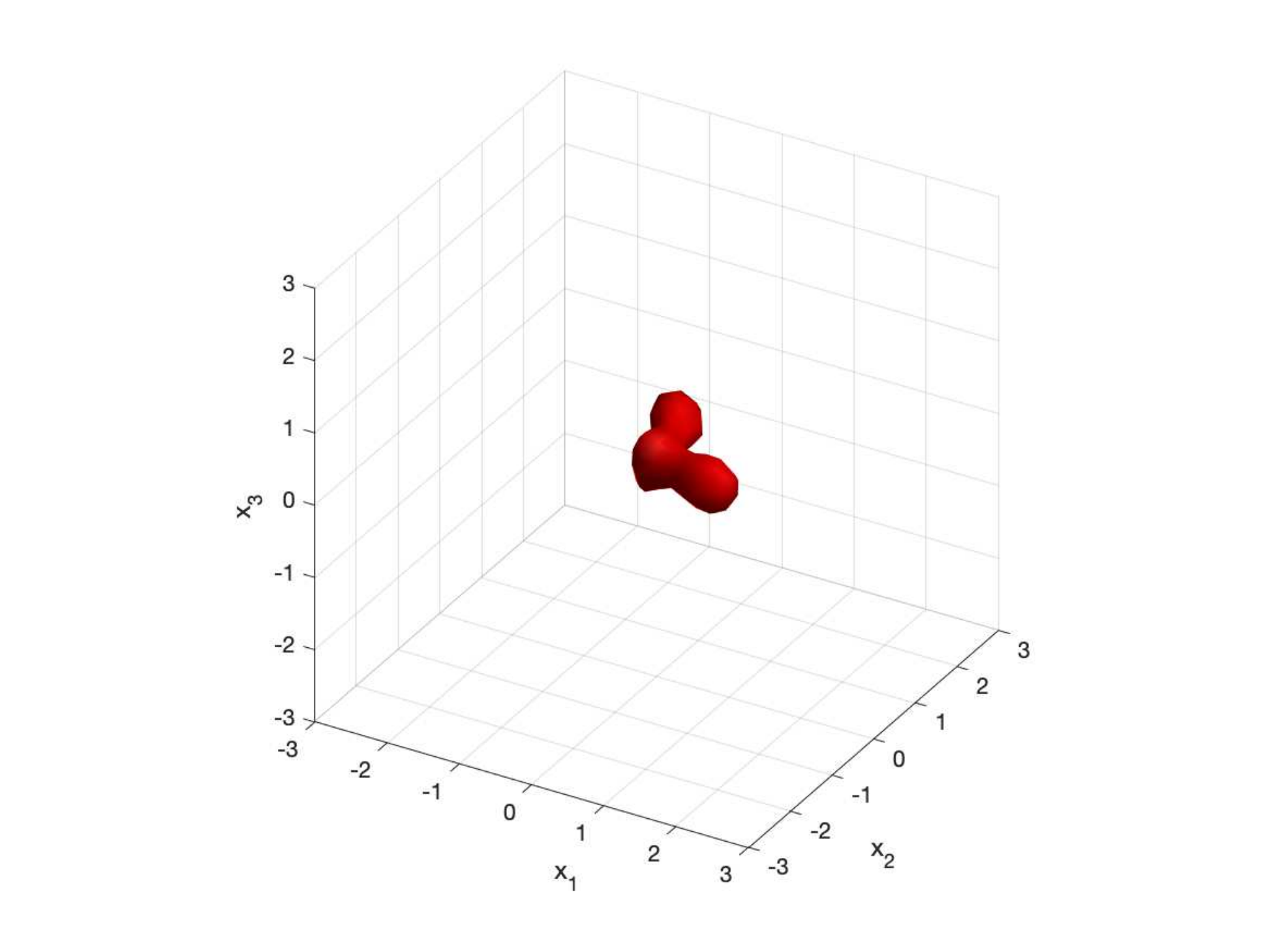}
    \includegraphics[width=5cm]{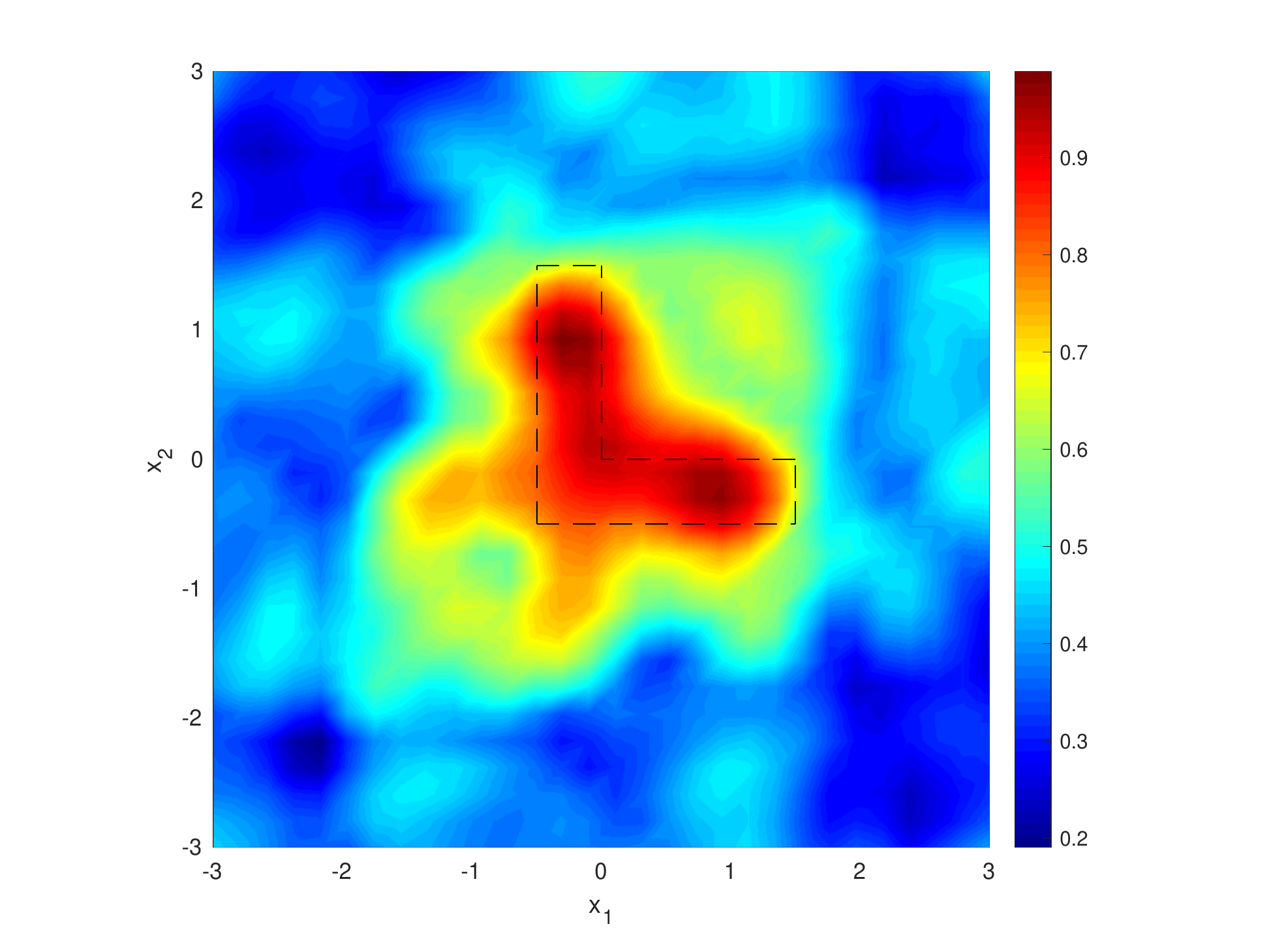}
      \includegraphics[width=5cm]{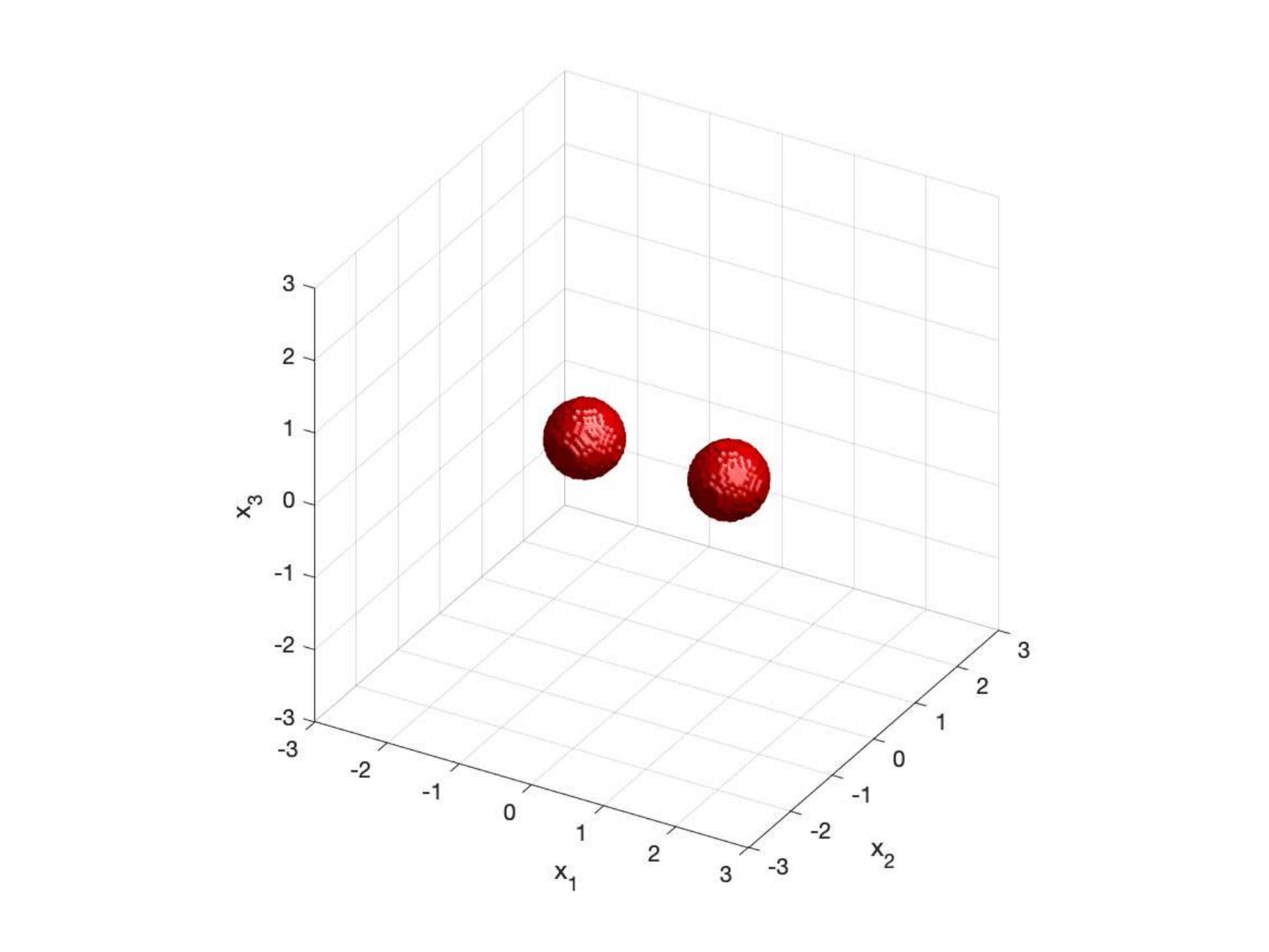}
    \includegraphics[width=5cm]{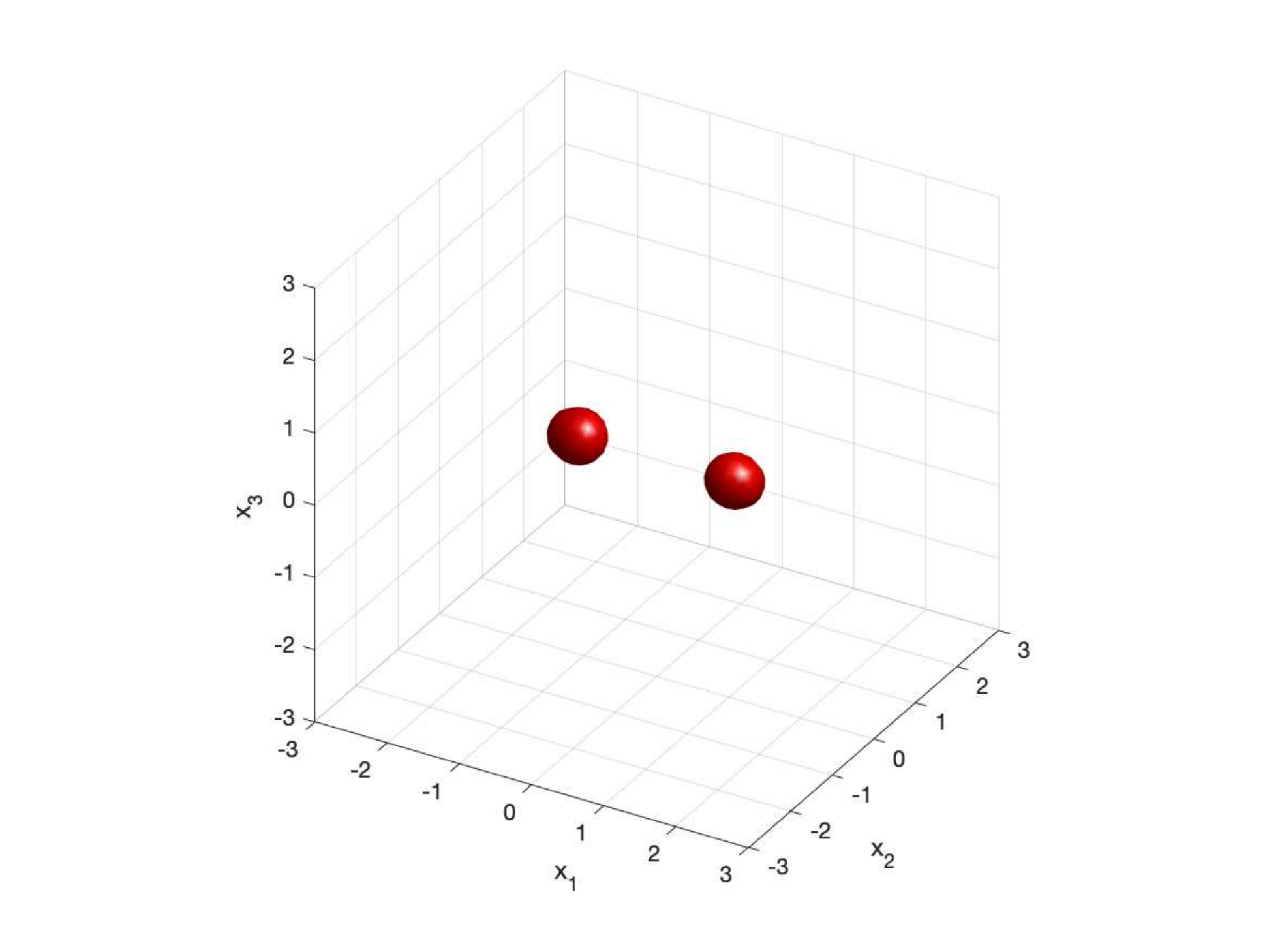}
    \includegraphics[width=5cm]{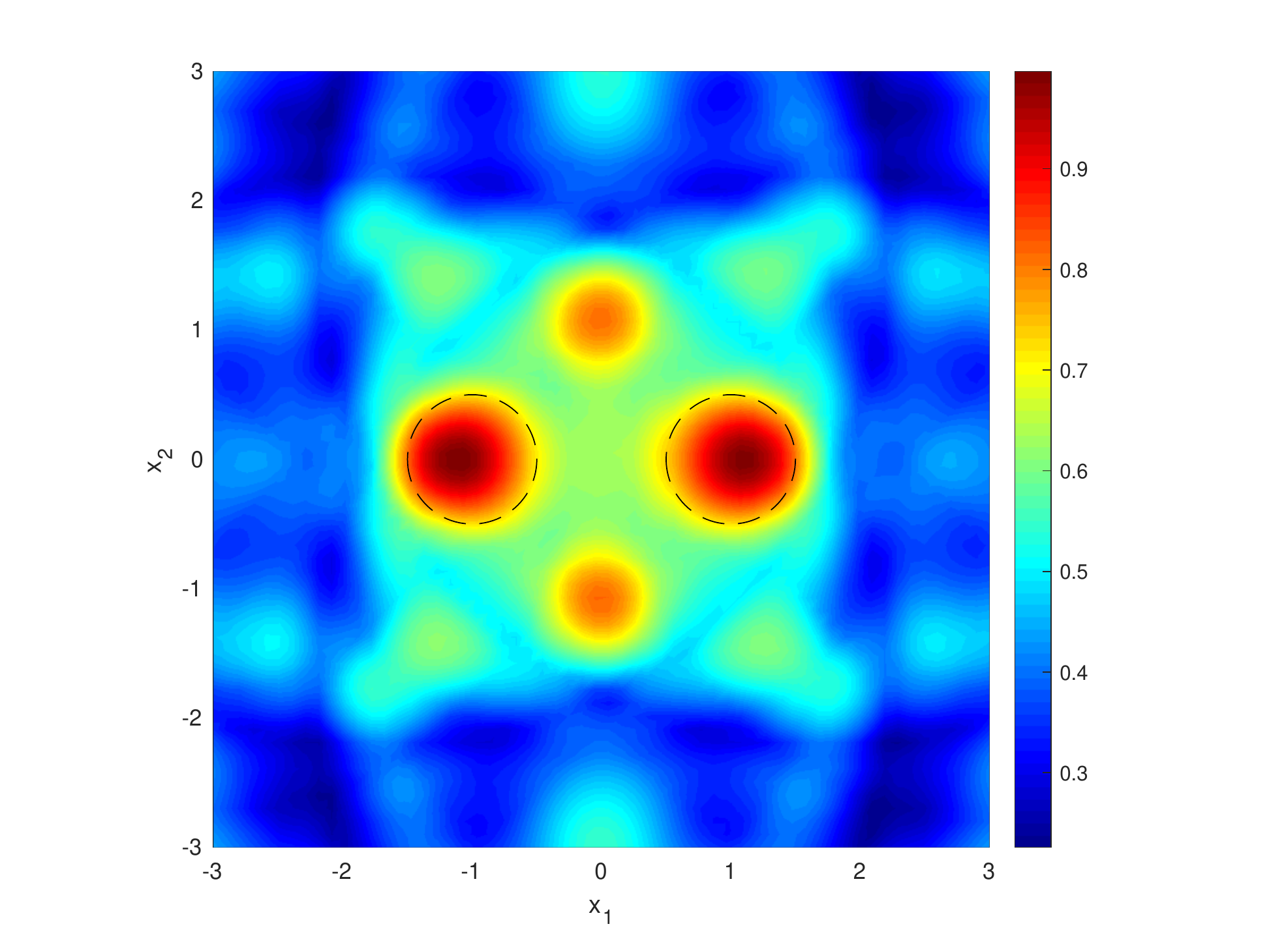}
\caption{{\bf Reconstructions of the L-shape (top) and two balls (bottom) with $14$ measurement points.} Left: exact support. Middle: iso-surface view of the reconstruction with iso-values $8.7\times 10^{-1}$  (L-shape) and $8.5\times 10^{-1}$ (two balls). Right: $x_1x_2$ cross section view of the reconstruction. }
\label{LBall pt 14}
\end{figure}

We continue to consider a variety of different geometries of the support: a cube (reconstruction in Figure \ref{cube pt 14}) given by
$$
\{ x=(x_1,x_2,x_3): |x_1| < 1, |x_2| < 1, |x_3| < 1 \},
$$
a rounded cylinder (reconstruction in Figure \ref{cr pt 14}) given by
$$
\Bigg\{ x=(x_1,x_2,x_3):
\begin{array}{ccc}
\sqrt{x_1^2 + x_2^2} <1  &  \mbox{ for } & |x_3| < 1  \\
\sqrt{x_1^2 + x_2^2 + (x_3-1)^2} <1  & \mbox{ for }  &  1 < x_3 < 2\\
\sqrt{x_1^2 + x_2^2 + (x_3+1)^2} <1  &  \mbox{ for }  &   -2  < x_3 <  -1
\end{array}\Bigg\},
$$
a peanut-shape support (reconstruction in Figure \ref{peanut pt 14}) given by
$$
\{ x=(x_1,x_2,x_3): \sqrt{(x_1-0.5)^2 + x_2^2+ x_3^2} < 1 \mbox{ or } \sqrt{(x_1+0.5)^2 + x_2^2+ x_3^2} < 1\},
$$
a L-shape support (reconstruction in Figure \ref{LBall pt 14}) given by
\begin{eqnarray*}
\{x: -0.5<x_1<0, -0.5<x_2<1.5, |x_3|<0.25\} \cup \{x: 0<x_1<1.5, -0.5<x_2<0, |x_3|<0.25\},
\end{eqnarray*}
and  two balls (reconstruction in Figure \ref{LBall pt 14}) given by
\begin{eqnarray*}
\{x: \sqrt{|x_1+1|^2 +  |x_2|^2 +  |x_3|^2}< 0.5 \} \cup\{x: \sqrt{|x_1-1|^2 +  |x_2|^2 +  |x_3|^2}< 0.5 \}.
\end{eqnarray*}


We observe from the numerical examples that (1) the annulus support of the source can be reconstructed from only one measurement point; (2) the location of a single source can be reconstructed from a few (sparse) measurement points (where we use $3$ measurement points with $11$ different frequencies); (3) both the shape and location of a single source can be well reconstructed from  multiple (but sparse) measurement points (where we use $14$ measurement points with $11$ different frequencies) that are roughly equally distributed all around the source. It is also observed that  the concave part of the support (see Figure \ref{peanut pt 14} for the peanut and top  of Figure \ref{LBall pt 14} for the L-shape) may be reconstructed; (4) the imaging algorithm also works when there are multiple sources (see bottom of Figure \ref{LBall pt 14}).

\section*{Acknowledgement}
The research of X. Liu is supported by the NNSF of China grant 11971471 and the Youth Innovation Promotion Association, CAS.

\section{Appendix: extension to far field measurements} \label{section appendix}
In this appendix we discuss briefly the indicator function defined by the multi-frequency sparse far field measurements. From the asymptotic behavior of Hankel functions, we deduce that the corresponding far field pattern of $u^s$ \eqref{us} has the form
\be\label{FarFieldrep}
u^{\infty}(\hx,k)=\int_{D}e^{-ik\hx\cdot\,y}f(y)dy,\quad \hx\in\,S^{2},\, k \in K.
\en
We consider the following multi-frequency  sparse far field measurements
\be\label{eqivalentdatafar}
\mathbb{M}_{F}:=\{u^{\infty}(\hx, k)\,|\, \hx\in\,\Theta_{L},\,\, k\in K\}
\en
with  $\Theta_{L}:=\{\pm\hx_1, \pm\hx_2, \cdots, \pm\hx_L\}\subset S^{2}$.

 In view of \eqref{FarFieldrep}, we have
\ben
u^{\infty}(\hx, -k) = u^{\infty}(-\hx, k), \quad \hx\in\Theta_L, \, k\in K.
\enn
Therefore, the multi-frequency sparse far field measurements $\mathbb{M}_{F}$ \eqref{eqivalentdatafar} gives the following data
\ben
u^{\infty}(\hx, k),\quad \hx\in\,\Theta_{L},\,\, k\in [-k_{max},k_{max}]\ba\{0\}.
\enn
We define the multi-frequency far field operator
$\mathcal {F}_{\hx}: L^2(K)\rightarrow L^2(K)$ by
\be\label{FarFieldOperator}
(\mathcal {F}_{\hx}\phi)(t) := \int_{K} u^{\infty}(\hx, t-s)\phi(s)ds, \quad t\in K.
\en
The analogous results of Theorem \ref{Thm-Nfac} and Theorem \ref{Itheorey} are formulated in the following theorem.
\begin{theorem}\label{Thm-Ffac}
The far field operator $\mathcal {F}_{\hx}: L^2(K)\rightarrow L^2(K)$ has a factorization in the form
\be\label{Ffac}
\mathcal {F}_{\hx}= \mathcal {Q}_{\hx}T_{\hx}\mathcal {Q}_{\hx}^{\ast}
\en
Here, $\mathcal {Q}_{\hx}: L^2(D)\rightarrow L^2(K)$ is given by
\ben
(\mathcal {Q}_{\hx}\psi)(t):= \int_{D}e^{-it\hx\cdot y}\psi(y)dy,\quad t\in K,
\enn
and its adjoint $\mathcal {Q}_{\hx}^{\ast}: L^2(K)\rightarrow L^2(D)$ is given by
\ben
(\mathcal {Q}_{\hx}^{\ast}\phi)(y):= \int_{K}e^{is\hx\cdot y}\phi(s)ds,\quad y\in D.
\enn
The operator $T_{\hx}: L^2(D)\rightarrow L^2(D)$ is a multiplication operator given by $T_{\hx}g=fg$, where  $f\in L^2(\R^3)$ is the  source with support $D$.
Moreover, it holds that
\ben
c_f\|\mathcal {Q}^{\ast}_{x}\phi\|^2_{L^{2}(D)} \leq |(\mathcal {F}_{\hx} \phi, \phi)_{L^2(K)}| \leq C_f\|\mathcal {Q}^{\ast}_{x}\phi\|^2_{L^{2}(D)}, \quad \forall\phi\in L^2(K),\, \hx\in\Theta_L,
\enn
where $c_f$ and $C_f$ are two constants given in \eqref{Assumpf}.
\end{theorem}
Motivated by the results in Theorem \ref{Thm-Ffac}, we introduce the following indicator function
\be\label{I-Far}
I(z):=\sum_{\hx\in\Theta_L}\Big|(\mathcal {F}_{\hx} \phi_{\hat{x}z}, \phi_{\hat{x}z})_{L^2(K)}\Big|,\quad z\in\R^3,
\en
where $\phi_{\hat{x}z}\in L^2(K)$ is given by
\be\label{phiz}
\phi_{\hat{x}z}(t)=e^{it\hx\cdot z}, \,t\in K.
\en
Consequently we have for $y\in D$ that
\be\label{ft}
\Big|(\mathcal {Q}^{\ast}_{x}\phi_{\hat{x}z})(y)\Big|=\Big|\int_{K}e^{is\tilde{\mathbbm{t}}}ds\Big|, \quad \tilde{\mathbbm{t}} = \hx\cdot(z-y),
\en
where $\phi_{\hat{x}z}$ is given by \eqref{phiz}. Similarly, $\Big|(\mathcal {Q}^{\ast}_{x}\phi_{\hat{x}z})(y)\Big|$ is a point spread function of $\tilde{\mathbbm{t}}$.

We remark that the factorization \eqref{Ffac} has been derived in \cite{GriesmaierSchmiedecke-source} where a factorization method is investigated. As can be seen, we have allowed a complex-valued  $f$ in the case of far field measurements. We finally remark that the results in this appendix can be directly extended to the two dimensional case.

\bibliographystyle{SIAM}

\end{document}